\documentclass[epjST]{svjour}

\usepackage{amsmath,amssymb}

\usepackage{amsthm}
\numberwithin{equation}{section}
\usepackage{graphicx}
\usepackage{geometry}
\usepackage{physics}
\usepackage{bm}
\usepackage{float}
\usepackage{xcolor}
\usepackage{enumitem}
\numberwithin{figure}{section}
\usepackage[most]{tcolorbox}
\usepackage{tikz}
\usetikzlibrary{positioning,arrows.meta,fit,calc}

\newtcolorbox{exbox}{
  colback=white,
  colframe=black,
  boxrule=0.5pt,
  left=4pt,
  right=4pt,
  top=4pt,
  bottom=4pt,
  arc=0pt,
  outer arc=0pt
}

\newcommand{\MatIn}[2]{\begingroup\raggedright\sloppy\noindent\texttt{In[#1]:= #2}\par\endgroup}
\newcommand{\MatOut}[2]{\begingroup\raggedright\sloppy\noindent\texttt{Out[#1]= #2}\par\endgroup}

\theoremstyle{plain}

\usepackage[
    backend=biber,
    style=numeric,
    sorting=none
]{biblatex}

\usepackage[
    colorlinks=true,
    linkcolor=blue!60!black,
citecolor=blue!60!black,
urlcolor=blue!60!black
]{hyperref}

\usepackage{listings}

\definecolor{mmaIn}{HTML}{4A6984}  
\definecolor{mmaOut}{HTML}{888888} 
\definecolor{mmaVar}{HTML}{0044B3} 

\newcommand{\fname}[1]{\subsection*{\hypertarget{cmd:#1}{#1}}\addcontentsline{toc}{subsection}{\texttt{#1}}}

\newcommand{\cmd}[1]{\hyperlink{cmd:#1}{\texttt{#1}}}
\newcommand{\cmdgroup}[2]{\noindent\emph{#1:} #2\par\vspace{3pt}}
\newcommand{\relcmds}[1]{%
\begin{tcolorbox}[colback=gray!5,colframe=black!60,boxrule=0.4pt,arc=0pt,
  left=5pt,right=5pt,top=4pt,bottom=2pt,breakable,before skip=12pt,after skip=12pt]
\small\raggedright
\textbf{Relevant \texttt{ECPack} commands} (each name links to its entry in Appendix~\ref{app:manual}).\par\smallskip
#1
\end{tcolorbox}}
\newcommand{\usageline}[1]{\begingroup\raggedright\sloppy\noindent\textbf{Usage: }\texttt{#1}\par\endgroup\smallskip}
\newcommand{\args}{\noindent\textbf{Arguments. }}
\newcommand{\returns}{\par\noindent\textbf{Returns.}}
\newcommand{\notes}{\par\noindent\textbf{Notes \& restrictions.}}
\newcommand{\slowmark}{$^{\dagger}$}

\lstdefinestyle{MathematicaStyle}{
    language=Mathematica,
    basicstyle=\ttfamily\small,
    numbers=none,                  
    frame=none,                    
    backgroundcolor=\color{white}, 
    xleftmargin=4em,               
    breaklines=true,
    breakatwhitespace=true,
    mathescape=true,               
    escapeinside={(*@}{@*)},       
    showspaces=false,
    showstringspaces=false
    columns=fullflexible,
    keepspaces=true,
    moredelim=[is][\color{mmaIn}]{|In|}{/In|},
    moredelim=[is][\color{mmaOut}]{|Out|}{/Out|}
}

\title{ECPack.wl: A Mathematica Package for Elliptic Curves over the Rationals and Finite Fields} 
\author{B. Ananthanarayan\inst{1}\fnmsep\thanks{\email{anant@iisc.ac.in}} \and 
Darsh R. Bhatt\inst{1}\fnmsep\thanks{\email{darshbhatt.ssa@gmail.com}}\thanks{Corresponding author.} \and 
Parastoo Taghizadeh Asl\inst{1}\fnmsep\thanks{\email{parastooasl@iisc.ac.in}}}

\institute{Centre for High Energy Physics, Indian Institute of Science, Bangalore 560 012, India}
\date{}

\begin{document}
\abstract{We present a Mathematica package entitled ECPack.wl, where we implement many important properties and commands relevant to the study of elliptic curves over the
rationals, and finite fields, namely $Q$ and $F_q$, where $q=p^r$ and $p$ is prime. By streamlining the analysis of general elliptic curves and their transformations, ECPack.wl offers foundational computational tools for applications in high-energy physics, including elliptic Feynman integrals.
This package significantly extends the functionalities of a package due to McGee in
so far as we are able to study all versions of elliptic curves such as the
long Weierstrass form, quartic, as well as Legendre forms.  We implement suitable
transformations that are required for the general case, as well as those associated
with the fields of characteristic 2 and 3.  Furthermore, the package computes many curve properties—such as point arithmetic, minimality checks, rank bounds and analytic rank, group structure, invariants, and the Weil and Tate pairings—over both Q and $F_p$, and the Kodaira classification of the singular fibres of one-parameter families of elliptic curves, with the associated F-theory data. 
The package and its accompanying notebooks are available at this URL: \href{https://github.com/DhamakedarDarsh/Mathematica-Code-and-worked-out-notebooks-for-Elliptic-Curves}{ECPack Github Repository}}

\maketitle
\begin{flushleft}

\textbf{Program Summary}\\
Program Title: ECPack.wl, version 1.0\\
Repository Link: https://github.com/DhamakedarDarsh/Mathematica-Code-and-worked-out-notebooks-for-Elliptic-Curves \\
Licensing Provisions: GNU General Public License\\
Programming Language: \texttt{Wolfram} Mathematica version 13.2 or higher

\end{flushleft}

\tableofcontents
\newpage
\section{Introduction}
In the recent past, elliptic curves, elliptic functions, elliptic integrals, and modular forms and their L-functions have begun to enter theoretical physics. For instance, in the calculation of certain Feynman integrals, elliptic and modular structures arise; hence, this category of Feynman integrals is often referred to as elliptic Feynman integrals. To calculate these integrals, we must study elliptic and modular structures. A recent collection of articles may be found in~\cite{BluemleinSchneiderPaule2019EllipticIntegralsQFT}, and 
for instance in Chapter 13 of the standard textbook~\cite{Weinzierl}. 
Furthermore, they arise in other areas of theoretical high-energy physics, such as string theory; see, e.g.~\cite{DHokerKaidi2024}.

In view of these applications, it is useful to study these structures systematically. In this article we focus on elliptic curves and introduce \texttt{ECPack.wl}, a Mathematica package for studying elliptic curves over the rationals and over finite fields. 
Note that other topics will be studied in our subsequent articles.
Although many physics applications ultimately lead to elliptic curves over \(\mathbb{C}\), here we confine ourselves to the rationals and finite fields as a first step. The transformations between cubic, Weierstrass, quartic, and Legendre forms are already useful in that broader setting, so this package is indirectly motivated by physics as well as by arithmetic geometry.

The subject of elliptic curves is one of the jewels of nineteenth-century mathematics, originating with Abel, Gauss, Jacobi, and Legendre as noted in the Introduction chapter of ref.~\cite{McKeanMoll1999}. These curves are associated with modular forms, topology, complex analysis, group theory, number theory, and modern cryptography. An elliptic curve over $\mathbb{C}$ is a smooth algebraic curve of genus one in $\mathbb{CP}^{2}$ together with a distinguished point, usually taken as the point at infinity. Elliptic curves can have points with coordinates in any field, such as $\mathbb{F}_{q}$, $\mathbb{Q}$, $\mathbb{R}$, or $\mathbb{C}$.\cite{Silverman2006IntroductionEllipticCurves} They are defined by an equation given by:
\begin{equation}
y^{2}+a_{1}xy+a_{3}y=x^{3}+a_{2}x^{2}+a_{4}x+a_{6}, \text{the long Weierstrass equation}
\end{equation}

Our purpose is to learn the use of elliptic curve technology from basic principles and to study various manifestations of these that is often found in standard textbooks.  Keeping in mind that this topic falls outside a typical university course syllabus, we wished to approach it through problem-solving and from the perspective of using computer algebra, especially the Mathematica platform.

To our knowledge, there is only one package by John McGee\cite{McGee2015}, which 
deals with elliptic curves and introduces several relevant topics, such as a linear combination of the two points, the order of the points, the group order, Hasse’s bound, the j-invariant of the curve, Bachet's duplication formula, as well as cryptography with ElGamal encoding and decoding using elliptic curve operators. The author has given a wishlist at the end of his accompanying notebooks, which lists Computation of the Weil Pairing $(P,Q)\mapsto e^{i\,2k\pi/n}$ and Implementation of the Schoof--Elkies--Atkin algorithm for $\left|E(\mathbb{F}_{q})\right|$, etc, to our knowledge, none of them have been made available to the public.  Miller's algorithm and the Weil and Tate pairings are implemented in \texttt{ECPack} (Section~\ref{sec:crypto}).  In this work, we begin by writing an independent package that, on the one hand, reproduces the computations demonstrated in the notebook of McGee~\cite{McGee2015} and, whenever possible, offers comparisons, and, on the other hand, significantly enhances it. Our aim in this work is to comprehensively study elliptic curves over rational and finite fields. Our package offers several extensions, such as 
\begin{enumerate}
    \item \textbf{Transformations and Coordinates}\\
Cubic to Weierstrass, quartic to Weierstrass, long to short Weierstrass, as well as Legendre to Weierstrass conversions, and coordinate conversions. \emph{The commands pertinent to transformations work for all elliptic curves regardless of the field over which they are defined.}

\item \textbf{Core Arithmetic and Equivalences}\\
Checking whether a point is on a curve or not, checking exact projective equivalence of two points, the chord and tangent operation, point addition, point multiplication by a scalar value, and point order.

\item \textbf{Minimal Models and Discriminants}\\
Computing the torsion subgroup of an elliptic curve, finding rational torsion points using the Nagell-Lutz theorem, finding independent generators for an elliptic curve, calculating a basic upper bound on a curve's rank, searching for points on a curve within a numerical bound, and Bachet's duplication.

\item \textbf{Heights and Pairings}\\
Computing the naive logarithmic height of a point, computing the canonical height using iterative limits, computing the height pairing of two points, computing the height matrix for a given list of points, and checking whether a set of points is independent based on their height matrix determinant or not.

\item \textbf{Finite Fields, ElGamal Cryptography and Pairings}\\
Adding two points modulo prime $p$, multiplying a point by scalar n modulo prime p, embedding an integer message onto the curve using Koblitz encoding, performing ElGamal encryption/encoding, and performing ElGamal decryption/decoding. Evaluating Miller's function, computing the Weil pairing and the reduced Tate pairing over $\mathbb{F}_p$ and $\mathbb{F}_{p^r}$, and finding a basis of the full $n$-torsion subgroup.

\item \textbf{Elliptic Families and Kodaira Classification}\\
Orders of vanishing of $c_4$, $c_6$ and $\Delta$ along a point of the base, the Kodaira type and Euler number of the singular fibre of a family of elliptic curves over $\mathbb{Q}(t)$, including the fibre at $t=\infty$, the list of all singular fibres with the check that the Euler numbers add up to $12n$, the split or non-split character of $I_n$ and $IV$ fibres, the $j$-function of a family, F-theory gauge algebras, scaling dimensions of Coulomb-branch parameters, Tate forms with prescribed vanishing orders, and the universal families with an $N$-torsion point.
\end{enumerate}

\begin{figure}[htbp]
\centering
\resizebox{\linewidth}{!}{%
\begin{tikzpicture}[
  font=\sffamily\footnotesize,
  >=Stealth,
  box/.style={draw, rounded corners=2pt, align=center, inner sep=4pt, text width=3.5cm, minimum height=0.9cm, fill=white},
  mcg/.style={box, fill=black!14},
  todo/.style={box, dashed},
  head/.style={draw, rounded corners=2pt, align=center, inner sep=5pt, fill=black!75, text=white, font=\sffamily\footnotesize\bfseries, text width=3.5cm, minimum height=1.1cm},
  wide/.style={box, text width=15.2cm},
  arr/.style={->, thick}
]
\node[wide] (inp) at (5.85,0) {\textbf{Input curve} (\texttt{curve}): cubic $\cdot$ long Weierstrass $\cdot$ short Weierstrass $\cdot$ quartic $\cdot$ Legendre $\cdot$ equation or coefficient list};
\node[wide, below=0.7cm of inp] (tr) {\textbf{Transformations and coordinates}: reduction to Weierstrass form, long $\leftrightarrow$ short, minimal isomorphisms, point-coordinate maps; separate routines for characteristic 2 and 3};
\draw[arr] (inp) -- (tr);

\node[head] (h1) at (0,-4.2) {Rational points\\$E(\mathbb{Q})$};
\node[head] (h2) at (3.9,-4.2) {Reduction and local data};
\node[head] (h3) at (7.8,-4.2) {Finite fields\\$\mathbb{F}_p,\ \mathbb{F}_{p^r}$};
\node[head] (h4) at (11.7,-4.2) {Cryptography};
\foreach \i in {h1,h2,h3,h4} \draw[arr] (\i |- tr.south) -- (\i.north);

\node[mcg,  below=3mm of h1] (a1) {Point addition, multiplication, order; Bachet duplication};
\node[box,  below=3mm of a1] (a2) {Torsion subgroup (Nagell--Lutz); bounded point search};
\node[box,  below=3mm of a2] (a3) {Generators; rank bounds; analytic rank};
\node[box,  below=3mm of a3] (a4) {Canonical height, height pairing and matrix; independence test};

\node[box,  below=3mm of h2] (b1) {Invariants $b_i,\ c_4,\ c_6,\ \Delta,\ j$};
\node[box,  below=3mm of b1] (b2) {Minimal models (Laska--Kraus); transport of points};
\node[box,  below=3mm of b2] (b3) {Tate's algorithm, conductor};
\node[box,  below=3mm of b3] (b4) {Good/bad reduction modulo $p$};
\node[box,  below=3mm of b4] (b5) {Kodaira fibres of families over $\mathbb{Q}(t)$; F-theory data};

\node[mcg,  below=3mm of h3] (c1) {Point addition and multiplication over $\mathbb{F}_p$};
\node[box,  below=3mm of c1] (c2) {Field arithmetic and points over $\mathbb{F}_{p^r}$};
\node[mcg,  below=3mm of c2] (c3) {Group order (Hasse bound); point order};
\node[box,  below=3mm of c3] (c4) {Group structure $E(\mathbb{F}_q)\cong\mathbb{Z}/n_1\times\mathbb{Z}/n_2$};

\node[box,  below=3mm of h4] (d1) {Koblitz embedding of messages};
\node[box,  below=3mm of d1] (d2) {Point counting using Schoof's algorithm};
\node[mcg,  below=3mm of d2] (d3) {ElGamal encryption and decryption over $\mathbb{F}_p$};
\node[box,  below=3mm of d3] (d4) {Miller function; Weil and Tate pairings; $n$-torsion basis};

\node[fit=(a4)(b5)(c4)(d4), inner sep=0pt] (all) {};
\node[mcg,  minimum height=3mm, minimum width=6mm, text width=6mm, anchor=north west] (l1) at ($(all.south west)+(0,-5mm)$) {};
\node[anchor=west] at (l1.east) {also available in McGee's package};
\node[box,  minimum height=3mm, minimum width=6mm, text width=6mm, anchor=north west] (l2) at ($(l1.north east)+(5.6cm,0)$) {};
\node[anchor=west] at (l2.east) {new in \texttt{ECPack}};

\end{tikzpicture}}
\caption{Overview of \texttt{ECPack.wl}. Any curve is first brought to Weierstrass form, then analyzed over $\mathbb{Q}$, locally (reduction modulo $p$), over finite fields, or used in cryptography. Shaded boxes indicate functionality that is also described for the package of McGee~\cite{McGee2015}; all other boxes are new.}
\label{fig:flowchart}
\end{figure}
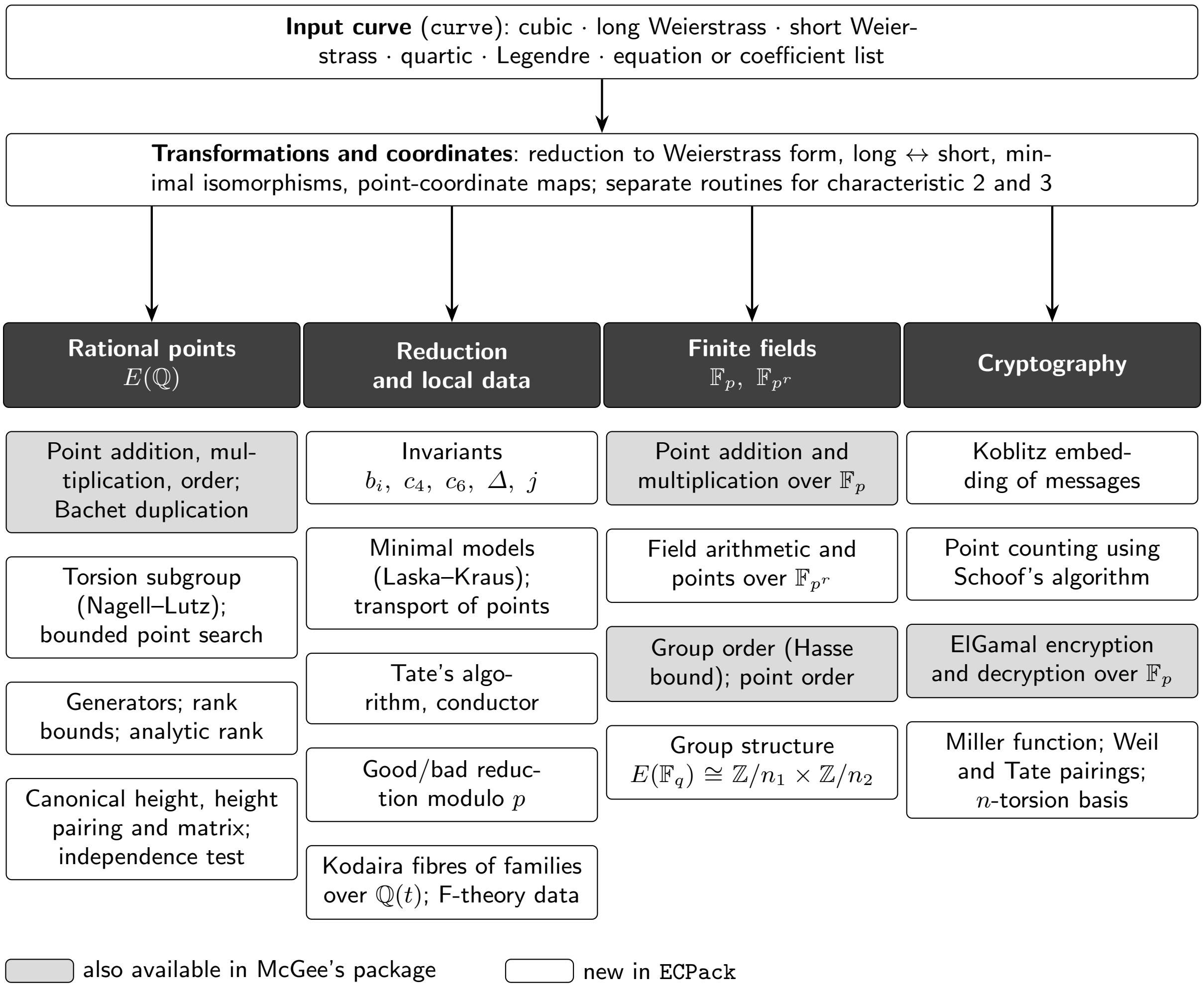

Figure~\ref{fig:flowchart} gives a schematic overview of the package: the functionality that is already available in the package of McGee~\cite{McGee2015} (as described above) is shaded, and everything else is new.

We are guided by the very useful and accessible book by Lozano-Robledo~\cite{LozanoRobledo2011} which allows us to give a solid structure to our code.

The first offering we have in this package is that we do not confine ourselves to the short Weierstrass form of the elliptic curve. For physical applications, one needs the quartic and Legendre forms, albeit over C. Nevertheless, our code now allows the study of elliptic curves in the following forms: cubic, long and short Weierstrass, quartic, and Legendre. We also code the requisite transformations, which are valid as long as one is not dealing with fields of characteristic 2 and 3.

The availability of good and bad reduction will also assist in evaluating the L-series associated with elliptic curves.

We also offer a realization of Schoof's point counting algorithm, which determines the trace of Frobenius modulo small primes and counts the points of a curve over a prime field $\mathbb{F}_p$ with $p$ far beyond the reach of enumeration.

At appropriate points, we refer to the accompanying notebooks, which collect textbook-style test cases from Lozano-Robledo, Husem{\"o}ller, Washington, Silverman, Silverman--Tate, McGee, and a set of general curves. These notebooks are meant to help readers see how the commands behave on representative inputs and to keep the main text focused on the ideas behind the package.

The rest of the article is organized as follows. Section~\ref{sec:forms} lists the forms in which an elliptic curve can be given and the maps to Weierstrass form. Section~\ref{sec:rationals} studies the group of rational points, its torsion, rank, and heights, and Section~\ref{sec:finite} studies the curve over finite fields, including point counting. Section~\ref{sec:reductions} treats reduction modulo $p$ and the minimal model, and Section~\ref{sec:kodaira} replaces the primes by the points of a parameter space and classifies the singular fibres of one-parameter families by Kodaira's table. Section~\ref{sec:crypto} turns to cryptography and pairings. The appendices recall the characteristic of a field (Appendix~\ref{app:fields}) and the Frobenius endomorphism (Appendix~\ref{app:frobenius}), and contain the package's user manual (Appendix~\ref{app:manual}) and notebook guide (Appendix~\ref{app: nb guide}).
We demonstrate that using our notebooks can significantly reduce tedium, and solutions can be readily obtained. Furthermore, we provide a list of all commands embedded in the code, which can be retrieved in the notebook using the \texttt{?} command.

\section{Various Forms of an Elliptic Curve}\label{sec:forms}

 In applications, an elliptic curve rarely arrives in the long Weierstrass form: curves that appear in the study of elliptic Feynman integrals can be given as a plane cubic or as a quartic, and families of curves are often written in the Legendre form. A first task of any package is therefore to recognize these forms and to bring each of them to a Weierstrass equation, on which all later computations are carried out. In this section we start by defining the elliptic curve and recall the forms accepted by \texttt{ECPack}, namely the general cubic, the long and short Weierstrass forms, the quartic and the Legendre form, and indicate how they are related. The transformations themselves are implemented by the commands listed at the end of the section.
 
\begin{definition}  

An elliptic curve over $\mathbb{C}$ is a smooth algebraic curve in $\mathbb{CP}^{2}$ of genus one together with one marked point, usually taken as the point at infinity.\normalfont\cite{Weinzierl}
The arithmetic of elliptic curves concerns itself with the points of a nonsingular cubic $y^2 = x^3 + ax + b$ with coordinates in a field, such as $\mathbb{F}_{p}$, $\mathbb{Q}$, $\mathbb{R}$, or $\mathbb{C}$.\normalfont\cite{McKeanMoll1999,Silverman2006IntroductionEllipticCurves}
\end{definition}

\begin{enumerate}

 \item\textbf{The General Cubic Form}

 An elliptic curve over $\mathbb{Q}$ is a smooth cubic projective curve $E$, with at least one rational point $\mathcal{O} \in E(\mathbb{Q})$. This point is the identity element of the group law defined over elliptic curves.
In the projective plane, it is given by a cubic polynomial $F(X,Y,Z)=0$ with rational coefficients \(a,b,c,\ldots \in \mathbb{Q}\), such that $E$ is smooth.
\begin{equation}
F(X,Y,Z)=aX^{3}+bX^{2}Y+cXY^{2}+dY^{3}+eX^{2}Z+fXYZ+gY^{2}Z+hXZ^{2}+jYZ^{2}+kZ^{3}=0.
\end{equation}

Regardless, we often consider just affine charts of E, e.g., those points of the form {[X, Y, 1]}, and study instead the affine curve given by, 

\begin{equation}
aX^{3}+bX^{2}Y+cXY^{2}+dY^{3}
+eX^{2}+fXY+gY^{2}
+hX+jY+k=0.
\end{equation}
This is the general cubic form of an elliptic curve.
Note that here, we left out some points of E at infinity (i.e., those points [X, Y, 0]).

\medskip

\item\textbf{The General Weierstrass Equation: The Long Weierstrass Form}

\begin{equation}
    y^{2}+a_{1}xy+a_{3}y=x^{3}+a_{2}x^{2}+a_{4}x+a_{6}
\end{equation}
With \(a_i \in \mathbb{Q}\).

\medskip

\item\textbf{The Short Weierstrass Form}

A curve E over a field K of characteristic (see Appendix~\ref{app:char}) different from 2 or 3, where $A,B\in K$, is given by a Weierstrass equation,

\begin{equation}
y^{2}=x^{3}+Ax+B.
\label{eq:short_weierstrass}
\end{equation}

This curve is non-singular if and only if

\begin{equation}
4A^{3}+27B^{2}\neq 0.
\label{eq:nonsingular}
\end{equation}

It also has a unique point at infinity $[0:1:0]$, which is called the \emph{origin} $\mathcal{O}$ or the \emph{point at infinity} of $E$.

\medskip

\item\textbf{Quartic Form}

We will meet elliptic curves defined by equations of the form
\begin{equation}
v^{2}=au^{4}+bu^{3}+cu^{2}+du+e,
\label{eq:quartic_general}
\end{equation}
with \(a\neq 0\). If we have a point \((p,q)\) lying on the curve with \(p,q\in K\), then the equation, when it is nonsingular, can be transformed into a Weierstrass equation by an invertible change of variables that uses rational functions with coefficients in the field \(K\). Note that an elliptic curve \(E\) defined over a field \(K\) always has a point in \(E(K)\), namely \(\infty\), whose projective coordinates \((0:1:0)\) certainly lie in \(K\). Therefore, if we are going to transform a curve \(C\) into Weierstrass form in such a way that all coefficients of the rational functions describing the transformation lie in \(K\), then we need to start with a point on \(C\) that has coordinates in \(K\).\cite{Washington2008EllipticCurves}\\

The quartic form can also be expressed as follows
    
    \begin{equation}
y^{2}
=
(x-x_{1})(x-x_{2})(x-x_{3})(x-x_{4}),
\label{eq:elliptic_curve_quartic_form}
\end{equation}

If all roots of the quartic polynomial on the right-hand side are distinct, this defines
a smooth elliptic curve.\cite{Weinzierl}

\medskip

\item\textbf{Legendre Form}

An elliptic curve in Legendre form over $K(\lambda)$ is defined by equations of the form
    \begin{equation}
y^{2}
=
x(x-1)(x-\lambda), 
\label{eq:elliptic_curve_legendre_form}
\end{equation} 
where \(K\) is a field of characteristic different from \(2\). The points \(\{(0,0),(1,0),(\lambda,0)\}\) are the three \(2\)-torsion points on \(E_{\lambda}\) for each value of \(\lambda\in K-\{0,1\}\).\\

Note that here \(\lambda\) is treated as an affine parameter, \(\lambda\in K\), so only \(\lambda=0,1\) are excluded. If we compactify the parameter space by adjoining the point at infinity, then \(\lambda\) is viewed as a projective parameter on \(\mathbb{P}^{1}(K)=K\cup\{\infty\}\). In that setting, the smooth Legendre family is parametrized by
\[
\lambda\in \mathbb{P}^{1}(K)-\{0,1,\infty\}.
\]
\end{enumerate}

\relcmds{
\cmdgroup{Reading and displaying a curve}{\cmd{ECPExtractCoefficients}, \cmd{ECPFormattedEquation}, \cmd{ECPToWeierstrassForm}, \cmd{ECPParseAffineCurve}, \cmd{ECPParseHomogeneousCurve}, \cmd{ECPPolynomialHomogeneous}, \cmd{ECPPolynomialTransform}, \cmd{ECPSetEC}}
\cmdgroup{Cubic, quartic and Legendre forms to Weierstrass form}{\cmd{ECPCubicToWeierstrass}, \cmd{ECPCubicToWeierstrassCharTwoThree}, \cmd{ECPQuarticToWeierstrass}, \cmd{ECPLegendreToWeierstrass}}
\cmdgroup{Long and short forms; normal forms in characteristic 2 and 3}{\cmd{ECPToShortWeierstrass}, \cmd{ECPShortWeierstrassForm}, \cmd{ECPToShortWeierstrassFp}, \cmd{ECPNormalFormCharTwoThree}}
\cmdgroup{Isomorphisms and point coordinates}{\cmd{ECPTransformWeierstrass}, \cmd{ECPComposeTransformations}, \cmd{ECPToShortCoordinates}, \cmd{ECPFromShortCoordinates}, \cmd{ECPToShortCoordinatesFp}, \cmd{ECPFromShortCoordinatesFp}, \cmd{ECPToLongCoordinates}}
\cmdgroup{Invariants}{\cmd{ECPWInvariants}, \cmd{ECPWeierstrassDiscriminant}, \cmd{ECPShortDiscriminant}, \cmd{ECPJInvariant}}
}

\section{Elliptic Curves over Rational Fields}\label{sec:rationals}

Once a curve has been brought to Weierstrass form by the transformations of Section~\ref{sec:forms}, the first field over which we study it is $\mathbb{Q}$. The set of rational points of an elliptic curve is not merely a set but an abelian group, and many commands of \texttt{ECPack} (point addition and multiplication, torsion points, generators, heights) implement aspects of this structure. In this section we first recall the group law on $E(\mathbb{Q})$, then describe the structure of the group, namely its torsion subgroup and its rank, and finally introduce the canonical height and the N\'eron--Tate pairing, which allow one to decide whether a given set of rational points is independent. We work with the short Weierstrass equation $y^2=x^3+Ax+B$, to which every curve over $\mathbb{Q}$ can be reduced.

\subsection{The Group Law on $E(\mathbb{Q})$}\label{The_Group_Structure_of_E(Q)}\label{Point_Addition_on_an_Elliptic_Curve}

Let $E: y^{2}=x^{3}+Ax+B$ with $A,B\in\mathbb{Q}$, and
\begin{equation}
E(\mathbb{Q})=\{(x,y)\in\mathbb{Q}^{2} : y^{2}=x^{3}+Ax+B\}\cup\{\mathcal{O}\},
\end{equation}
where $\mathcal{O}=[0:1:0]$ is the point at infinity. $E(\mathbb{Q})$ carries a group structure whose operation is the addition of points described below.

\begin{figure}[H]
    \centering
    \begin{minipage}[t]{0.325\linewidth}
        \centering
        \includegraphics[width=\linewidth]{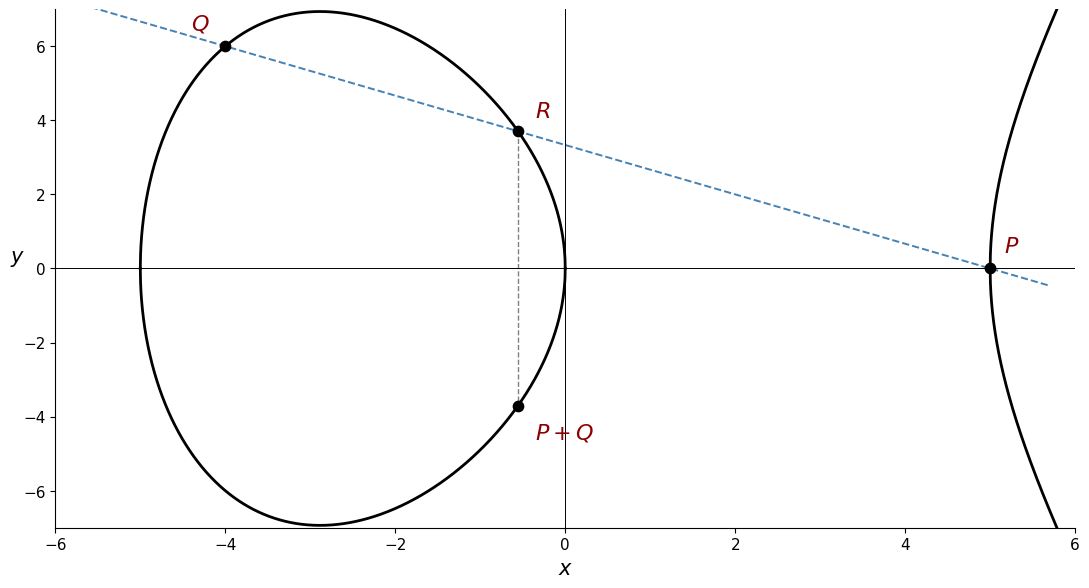}\\[2pt]
        (a) $P\neq Q$
    \end{minipage}\hfill
    \begin{minipage}[t]{0.325\linewidth}
        \centering
        \includegraphics[width=\linewidth]{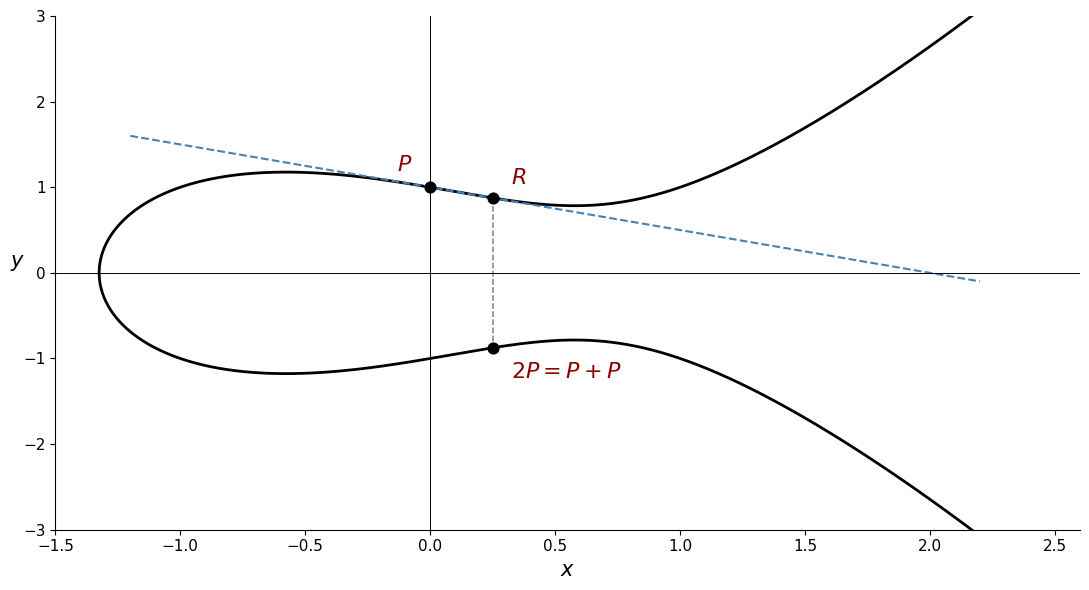}\\[2pt]
        (b) $P=Q$
    \end{minipage}\hfill
    \begin{minipage}[t]{0.325\linewidth}
        \centering
        \includegraphics[width=\linewidth]{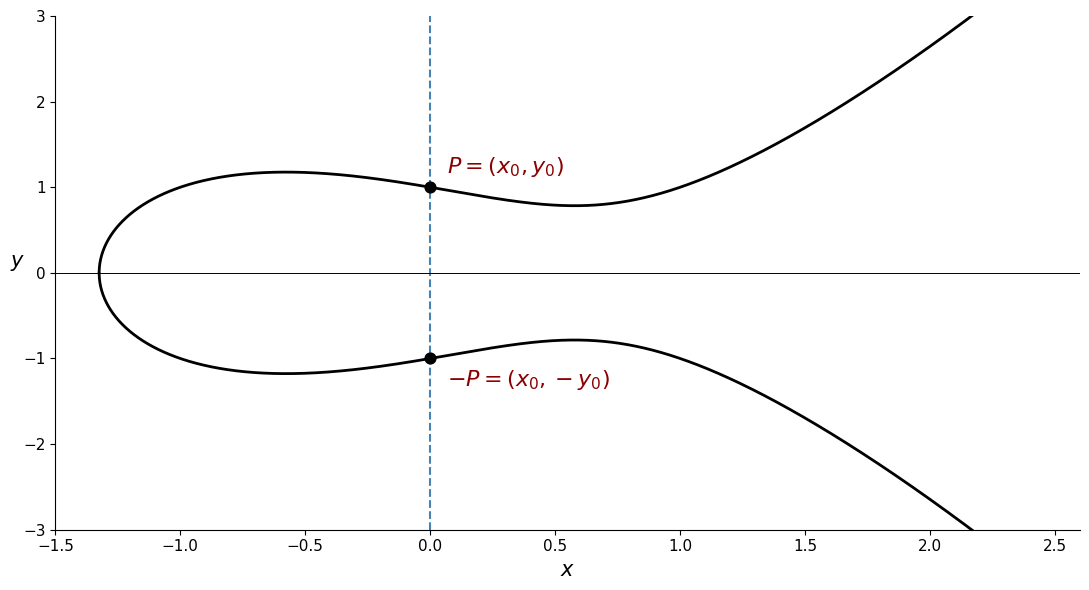}\\[2pt]
        (c) $P+(-P)=\mathcal{O}$
    \end{minipage}
    \caption{The group law on an elliptic curve: (a) addition of two distinct points, (b) adding a point to itself, (c) the vertical line through $P$ and $-P$, which explains the role of the point at infinity.}
    \label{fig:group_law}
\end{figure}

\paragraph{Chord and tangent.}
For two distinct points $P,Q\in E(\mathbb{Q})$, the line through $P$ and $Q$ meets $E$ in a third point $R$; the sum $P+Q$ is the reflection of $R$ across the $x$-axis (Fig.~\ref{fig:group_law}a). If $P=Q$, the line is replaced by the tangent to $E$ at $P$, which meets the curve at a new point whose reflection is $2P$ (Fig.~\ref{fig:group_law}b). The vertical line through $P$ and $-P$ meets $E$ in no third affine point. To define $P+(-P)$ we therefore adjoin the point $\mathcal{O}$ ``at infinity'', which lies on every vertical line, so that $P+(-P)=\mathcal{O}$ (Fig.~\ref{fig:group_law}c).\cite{Silverman2006IntroductionEllipticCurves}

\paragraph{Formulas.}
Let $P=(x_{1},y_{1})$ and $Q=(x_{2},y_{2})$, and let the line through them be $L: y=\lambda x+\nu$, with
\begin{equation}
    \lambda=
\begin{cases}
\dfrac{y_{2}-y_{1}}{x_{2}-x_{1}}, & \text{if } P\neq Q, \\[10pt]
\dfrac{3x_{1}^{2}+A}{2y_{1}}, & \text{if } P=Q,
\end{cases}
\qquad
\nu=y_{1}-\lambda x_{1}.
\end{equation}
Then $P+Q=(x_{3},y_{3})$ with
\begin{equation}
    x_{3}=\lambda^{2}-x_{1}-x_{2},\qquad y_{3}=-(\lambda x_{3}+\nu),
\end{equation}
and $-(x,y)=(x,-y)$. A vertical chord ($x_{1}=x_{2}$, $y_{1}=-y_{2}$) gives $\mathcal{O}$.

\begin{theorem}\normalfont\cite{Silverman2006IntroductionEllipticCurves}\label{thm:group_law}
With this addition, $E(\mathbb{Q})$ is an abelian group with identity $\mathcal{O}$: it is closed under addition, the operation is associative and commutative, $P+\mathcal{O}=P$ and $P+(-P)=\mathcal{O}$ for all $P\in E(\mathbb{Q})$.
\end{theorem}

\paragraph{Bachet's duplication formula.}
For $P=Q$ the formulas above reduce to Bachet's duplication formula (1621), found while expressing an integer as the difference of a square and a cube in rational numbers,
\begin{equation}
    y^2-x^3=c, \qquad c\in\mathbb{Z}.
\end{equation}
If $(x,y)\in \mathbb{Q}^2$ with $y\ne0$ is a solution, then
\begin{equation}
  \left(
\frac{x^{4}-8cx}{4y^{2}},
\frac{-x^{6}-20cx^{3}+8c^{2}}{8y^{3}}
\right)
\end{equation}
is also a rational solution. Hence if $c\notin\{1,-432\}$ and $xy\ne0$, the integer $c$ can be written as a difference of a square and a cube of nonzero rational numbers in infinitely many ways.\cite{silverman}

\subsection{Structure of $E(\mathbb{Q})$: Torsion and Rank}\label{subsec:structure}\label{subsec:primebound}

\begin{theorem}[\(E(\mathbb{Q})\): Mordell--Weil Group]\normalfont\cite{LozanoRobledo2011}\label{Mordell_Weil_Group}
\(E(\mathbb{Q})\) is a finitely generated abelian group. In other words, there are points
\(P_{1},\ldots,P_{n}\) such that any other point \(Q\in E(\mathbb{Q})\) can be expressed as a linear combination
\[Q=a_{1}P_{1}+a_{2}P_{2}+\cdots+a_{n}P_{n},
\qquad \text{for some } a_i\in\mathbb{Z}.\]
\end{theorem}

\begin{figure}[H]
    \centering
    \includegraphics[width=0.6\linewidth]{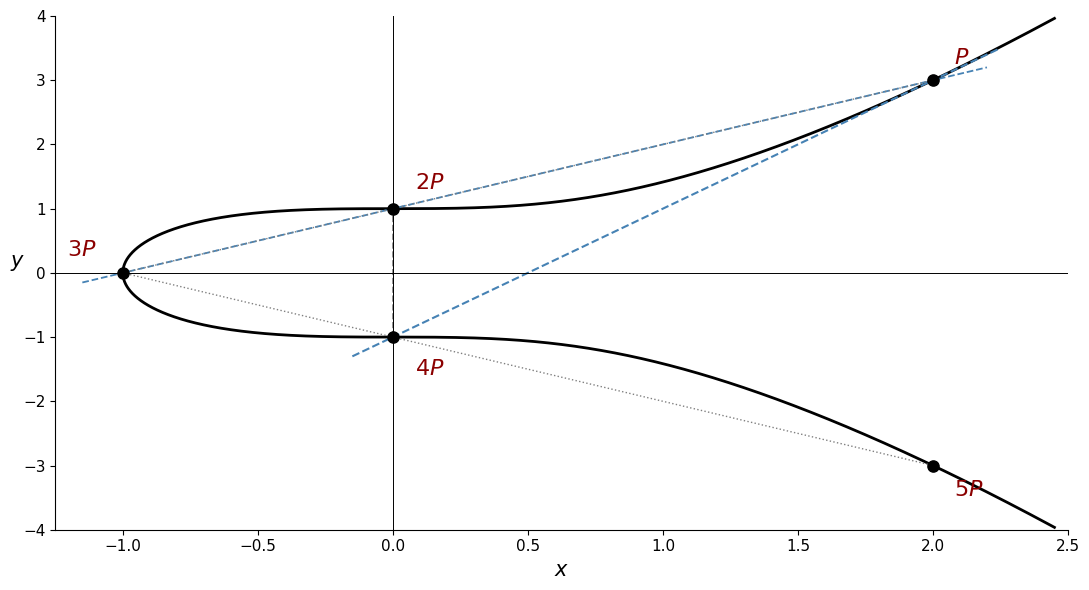}
    \caption{The rational torsion points of $y^2=x^3+1$, where $P=(2,3)$ generates $E(\mathbb{Q})_{\mathrm{torsion}}\cong\mathbb{Z}/6\mathbb{Z}$; the point $6P=\mathcal{O}$ at infinity is not shown.}
    \label{fig:The_rational_torsion_points}
\end{figure}

By the Mordell--Weil theorem and the structure theory of finitely generated abelian groups,
\begin{equation}
E(\mathbb{Q}) \cong E(\mathbb{Q})_{\mathrm{torsion}} \oplus \mathbb{Z}^{R_E},
\qquad
E(\mathbb{Q})_{\mathrm{torsion}}
=
\left\{
P\in E(\mathbb{Q})
\;\middle|\;
nP=\mathcal{O}\text{ for some } n\in\mathbb{N}
\right\}.
\end{equation}
The decomposition is not canonical. The torsion subgroup is the finite group of points of finite order. The free part $\mathbb{Z}^{R_E}$, with $R_E\geq 0$, is generated by $R_E$ independent points of infinite order, that is, points $P$ with $nP\neq \mathcal{O}$ for all nonzero $n\in \mathbb{Z}$. The integer $R_E$ is the \emph{rank} of $E/\mathbb{Q}$. The rational torsion points of an example curve are shown in Fig.~\ref{fig:The_rational_torsion_points}.

\begin{theorem}[Mazur]\normalfont\cite{LozanoRobledo2011}
For an elliptic curve \(E/\mathbb{Q}\), the
\(E(\mathbb{Q})_{\mathrm{torsion}}\) is isomorphic to exactly one of the following groups:
\[
\begin{array}{lll}
\mathbb{Z}/N\mathbb{Z}
& \text{with}
& 1\leq N\leq 10 \text{ or } N=12,\\[0.3cm]
\mathbb{Z}/2\mathbb{Z}\oplus\mathbb{Z}/2M\mathbb{Z}
& \text{with}
& 1\leq M\leq 4.
\end{array}
\]
\end{theorem}

In particular, if the order of a rational point \(P\in E(\mathbb{Q})\) is greater than \(12\), then \(P\) has infinite order, and \(E(\mathbb{Q})\) contains infinitely many distinct rational points.

Although the rank \(R_E\) cannot be read off from the equation, it admits upper bounds when \(E\) has a rational point of order \(2\). Let \(\nu(n)\) be the number of distinct primes dividing a nonzero integer \(n\), so that \(\nu(360)=3\) and \(\nu(\pm1)=0\).

\begin{theorem}[Rank bound]\normalfont\cite{LozanoRobledo2011}\label{thm:rankbound}
Let \(E:\ y^2=x^3+ax^2+bx\) with \(a,b\in\mathbb{Z}\). Then
\begin{equation}
R_E\;\le\;\nu(a^2-4b)+\nu(b)-1 .
\end{equation}
More generally, if \(E/\mathbb{Q}\) has a nontrivial rational \(2\)-torsion point, then \(R_E\le m+2\alpha-1\), where \(m\) and \(\alpha\) are the numbers of primes at which \(E\) has multiplicative and additive bad reduction, respectively.
\end{theorem}

\begin{proposition}\normalfont\cite{LozanoRobledo2011}\label{prop:threeroots}
If \(E:\ y^2=(x-e_1)(x-e_2)(x-e_3)\) with \(e_1,e_2,e_3\in\mathbb{Z}\), then \(R_E\le 2\,\nu(\Delta_E)\).
\end{proposition}

A rational \(2\)-torsion point exists exactly when the cubic has a rational root, and moving that root to \(x=0\) puts \(E\) in the form of Theorem~\ref{thm:rankbound}. The command \texttt{ECPRankBound} evaluates all bounds that apply and returns the smallest. These are upper bounds only: for \(y^2=x^3-x\) the bound is \(0\) and the rank is \(0\), whereas for \(y^2=x^3-1156x\) the bound is \(3\) and the rank is \(2\)~\cite{LozanoRobledo2011}. For a curve without a rational point of order \(2\), such as 37a1, the command returns \texttt{Missing["NoRational2Torsion"]}. The command \texttt{ECPRankUpperBound} is only a heuristic and not a theorem (it overestimates the rank of \(y^2=x^3-x\)).

\subsection{Heights and Linear Independence of Rational Points}\label{subsec:heights}
The material follows \cite{LozanoRobledo2011}.\\
Rational points \(P_1,\dots,P_m \in E(\mathbb{Q})\) are \emph{linearly independent} if
\begin{equation}
    n_1P_1+n_2P_2+\dots+n_mP_m=T, 
\end{equation}
with \(n_i \in \mathbb{Z}\) and \(T\) a torsion point, forces \(n_1=\dots=n_m=0\). To test this, we use heights. The height of \(\frac{m}{n}\in\mathbb{Q}\), with \(\gcd(m,n)=1\), is
\(h\left(\frac{m}{n}\right)=\log\bigl(\max\{|m|,|n|\}\bigr)\), and for a point \(P=(x,y)\in E(\mathbb{Q})\) we set \(H(P)=h(x)\). The \textbf{canonical height} of \(P\in E(\mathbb{Q})\) is
\begin{equation}
    \widehat{h}(P)
=
\frac{1}{2}
\lim_{N\to\infty}
\frac{H(2^{N}P)}{4^{N}},
\end{equation}
where \(2^N P\) is computed with the group law. For \(P,Q\in E(\mathbb{Q})\) and \(m\in\mathbb{Z}\) it satisfies
\begin{equation}
    \widehat{h}(P+Q)+\widehat{h}(P-Q)=2\widehat{h}(P)+2\widehat{h}(Q),
    \qquad
    \widehat{h}(mP)=m^{2}\widehat{h}(P),
    \qquad
    \widehat{h}(P)\geq 0,
\end{equation}
with \(\widehat{h}(P)=0\) if and only if \(P\) is a torsion point.

The \textbf{N\'eron--Tate pairing} is defined by
\begin{equation}
    \langle P,Q\rangle
=
\tfrac{1}{2}\bigl(\widehat{h}(P+Q)-\widehat{h}(P)-\widehat{h}(Q)\bigr).
\end{equation}
It is symmetric and bilinear, $\langle P,P\rangle=\widehat{h}(P)$, and so $\langle P,P\rangle=0$ exactly for torsion points: this distinguishes torsion points from points of infinite order. (This is the normalisation used by \texttt{ECPHeightPairing}, see Appendix~\ref{app:manual}.) For points \(P_1,\dots,P_r\) the \textbf{elliptic height matrix} is
\begin{equation}
    \mathcal{H}
=
\left(
\langle P_i,P_j\rangle
\right)_{1\leq i,j\leq r}.
\end{equation}

\begin{theorem}\normalfont\cite{LozanoRobledo2011}\label{thm:independence}
\begin{enumerate}
    \item[(a)] If \(\det(\mathcal{H})\neq 0\), the points \(\{P_i\}_{i=1}^{r}\) are linearly independent and the rank of \(E(\mathbb{Q})\) is at least \(r\).
    \item[(b)] If \(\det(\mathcal{H})=0\) and \(0\neq(n_1,\ldots,n_r)\in\mathbb{Z}^r\) lies in \(\ker(\mathcal{H})\), then the points are linearly dependent and \(\sum_{k=1}^{r} n_k P_k = T\) for a torsion point \(T\in E(\mathbb{Q})\).
\end{enumerate}
\end{theorem}

\relcmds{
\cmdgroup{Point arithmetic}{\cmd{ECPAdd}, \cmd{ECPMultiply}, \cmd{ECPChordAndTangent}, \cmd{ECPBachetDuplication}, \cmd{ECPOnCurveQ}, \cmd{ECPPointQ}, \cmd{ECPAreProjectivelyEquivalent}}
\cmdgroup{Torsion and point orders}{\cmd{ECPFindTorsionPoints}, \cmd{ECPTorsionSubgroup}, \cmd{ECPPointOrder}, \cmd{ECPFindPointOrder}}
\cmdgroup{Searching for points, generators and rank}{\cmd{ECPSearchPoints}, \cmd{ECPFindGenerators}, \cmd{ECPAnalyticRank}}
\cmdgroup{Heights, pairings and independence}{\cmd{ECPNaiveHeight}, \cmd{ECPCanonicalHeight}, \cmd{ECPHeightPairing}, \cmd{ECPHeightMatrix}, \cmd{ECPIndependentQ}}
\cmdgroup{Prime divisors and rank bounds}{\cmd{ECPPrimeDivisorCount}, \cmd{ECPRankBound}, \cmd{ECPLocalData}, \cmd{ECPRankUpperBound}}
}

\section{Elliptic Curves over Finite Fields}\label{sec:finite}

Elliptic curves over finite fields are one of the places where the geometry of the curve and the arithmetic of its points meet in a concrete way. Over $\mathbb{Q}$, one often studies individual rational points and the group law on the curve. Over a finite field $\mathbb{F}_q$, the picture changes: the set of points is finite, and an important question becomes not whether the points exist, but how many there are. This counting problem plays a central role in algorithmic number theory and elliptic curve cryptography.
\subsection{Curve structure}
An elliptic curve over \(\mathbb{F}_{q}\), where \(q=p^{r}\) for a prime $p$, in short Weierstrass form\cite{Washington2008EllipticCurves} is as follows
\begin{equation}
    E:\qquad y^{2}=x^{3}+Ax+B,
\qquad
A,B\in \mathbb{F}_{q},
\end{equation}
For a nonsingular curve\cite{Washington2008EllipticCurves}
\begin{equation}
    4A^{3}+27B^{2}\neq 0.
\end{equation}
Equivalently, the discriminant\cite{Washington2008EllipticCurves}
\begin{equation}
    \Delta=-16(4A^{3}+27B^{2})
\end{equation}
is nonzero as well.

\(E(\mathbb{F}_{q})\) is the set of all affine points of \(E\) whose coordinates belong to \(\mathbb{F}_{q}\), together with the point at infinity:
\begin{equation}
    E(\mathbb{F}_{q})
=
\{(x,y)\in \mathbb{F}_{q}^{2}: y^{2}=x^{3}+Ax+B\}\cup\{\mathcal{O}\}.
\end{equation}

Note that
\begin{itemize}
    \item When a curve is defined over \(\mathbb{F}_{q}\), it means that its defining equation has coefficients $A,B$ in \(\mathbb{F}_{q}\).
    \item A point is rational over \(\mathbb{F}_{q}\) if its coordinates lie in \(\mathbb{F}_{q}\).
    \item A curve defined over \(\mathbb{F}_{q}\) can also have points whose coordinates lie in a larger field. Thus, for every \(n\geq 1\) (see Appendix~\ref{app:prime}),
\begin{equation}
    E(\mathbb{F}_{q})\subseteq E(\mathbb{F}_{q^{n}}).
\end{equation}
\end{itemize}

\subsection{Advanced Point Counting}
A simple heuristic suggests that \(N_q=\#E(\mathbb{F}_{q})\) the number of $\mathbb{F}_{q}$-rational points on $E$ should be close to $q+1$, since for each $x\in\mathbb{F}_{q}$ there are at most two possible values of $y$, and the point at infinity contributes one extra point.

\begin{definition}[The Frobenius Endomorphism]\normalfont\cite{Silverman2006IntroductionEllipticCurves} (see Appendix~\ref{app:frobenius})\\
    Suppose $E$ is defined over the finite field $\mathbb{F}_q$:
\begin{equation}
    \phi_q:E\to E,\qquad \phi_q(x,y)=(x^q,y^q)
\end{equation}
is the natural map associated with the finite-field structure, and it encodes the curve's arithmetic very efficiently. From this point of view, \textbf{point counting is closely related to understanding how Frobenius acts on the curve.} This is the reason the trace $a_q$ appears so naturally in both theory and computation.
\end{definition}{}

\textbf{The Trace of Frobenius}\normalfont~\cite{Silverman2006IntroductionEllipticCurves} (see Appendix~\ref{app:trace}): The quantity \begin{equation}
    a_q = q + 1 - \#E(\mathbb{F}_q)
\end{equation} is called the trace of Frobenius, because one way to calculate it is to use the Frobenius map to obtain a linear transformation on a certain vector space \(V_{\ell}(E)\) that is associated with the Tate module. Then a$_q$ is the trace of that linear transformation.

\begin{theorem}[Hasse's theorem]\normalfont\cite{Washington2008EllipticCurves}\label{Theorem: Hasses's_Theorem} makes this precise:
\[
\bigl|N_q-(q+1)\bigr|\le 2\sqrt{q},
\]
So, although the exact count depends on the curve, it is always constrained to lie in a fairly small interval around $q+1$.
\end{theorem}

It is often convenient to package the deviation from the expected value into the trace of Frobenius,
\[
a_q=q+1-N_q.
\]
Then Hasse's bound becomes simply $|a_q|\le 2\sqrt{q}$.\\

For small fields, one can count points directly by checking each $x\in\mathbb{F}_q$ and determining whether the right-hand side of the curve equation is a quadratic residue. This works well for small examples and for verifying formulas by hand, \textbf{but it quickly becomes too slow when $q$ is large.} This is where advanced point-counting algorithms come into play.

\subsection{Schoof's Algorithm}\label{sec:schoof}
Schoof's algorithm \cite{Schoof1985} provides a breakthrough technique for computing \(N_q\). Using this technique, we need not directly compute the trace of Frobenius \(a_q\); instead, we compute the trace of Frobenius modulo many small primes \(\ell\ne char(\mathbb{F}_q)\) and then recover $a_q$ via the Chinese Remainder Theorem.\\
So, the problem becomes one of studying how Frobenius acts on the $\ell$-torsion subgroup $E[\ell]$. Because on the \(\ell\)-torsion points,
\(
E[\ell],
\)
Frobenius satisfies
\(
\phi_q^{2}-a_q\phi_q+q=0.
\) Modulo \(\ell\), this becomes a finite computation.

For each small prime \(\ell\neq char(\mathbb{F}_q)\), Schoof determines
\[
a_q \pmod{\ell}.
\]
The congruences for \(a_q\), obtained modulo different primes \(\ell\), are combined using the Chinese Remainder Theorem. Since the moduli are pairwise coprime, the combined modulus is the product of the individual moduli. Once this product of small primes \(\ell\) is larger than the width of Hasse's interval \(4\sqrt{q}\), Hasse's bound
\(
|a_q| \leq 2\sqrt{q}
\)
determines a unique value of \(a_q\), since there can be only one integer in the Hasse interval satisfying all the congruences. Finally,
\begin{equation}
    \#E(\mathbb{F}_q)=q+1-a_q.
\end{equation}

\paragraph{The computation modulo $\ell$.} Let $q=p$ be prime, $E: y^2=x^3+Ax+B$ and $\ell\neq p$ an odd prime. The $x$-coordinates of the nonzero points of $E[\ell]$ are the roots of the division polynomial $\psi_\ell\in\mathbb{F}_p[x]$, of degree $(\ell^2-1)/2$. All computations are carried out in the ring $\mathbb{F}_p[x]/(\psi_\ell)$, in which $P=(x,y)$ is a generic point of $E[\ell]$ and the Frobenius acts as $\phi(P)=\bigl(x^{p},\,(x^3+Ax+B)^{(p-1)/2}\,y\bigr)$. With $\bar q=p \bmod \ell$, one adds $\phi^2(P)$ and $[\bar q]P$ with the chord-and-tangent formulas of Section~\ref{sec:rationals} and looks for $\tau\in\{1,\dots,\ell-1\}$ with
\begin{equation}
    \phi^{2}(P)+[\bar q]P=[\tau]\,\phi(P);
\end{equation}
then $a_p\equiv\tau \pmod{\ell}$. If $\phi^2(P)$ and $[\bar q]P$ have the same $x$-coordinate for some $P\neq\mathcal{O}$, then $\phi^2(P)=\pm[\bar q]P$ and there are two cases: either $a_p\equiv0$, or $\bar q\equiv w^2 \pmod{\ell}$ is a square, $\phi(P)=\pm[w]P$ and $a_p\equiv\pm2w \pmod{\ell}$, the sign being read off from the $y$-coordinates. For $\ell=2$ the answer is simpler: $a_p$ is even if and only if $x^3+Ax+B$ has a root in $\mathbb{F}_p$, that is, if and only if $\gcd(x^p-x,\,x^3+Ax+B)\neq1$. All powers $x^p$ are computed by repeated squaring modulo $\psi_\ell$, and the primes $\ell$ needed to exceed $4\sqrt{p}$ are of size $O(\log p)$, so the running time is polynomial in $\log p$.

\paragraph{Implementation in \texttt{ECPack}.} \texttt{ECPSchoofPointCount} implements the algorithm above for a curve over a prime field $\mathbb{F}_p$ in any of the accepted forms, with arithmetic modulo $\psi_\ell$ done by the polynomial routines of the Wolfram Language over $\mathbb{Z}/p\mathbb{Z}$. It uses the primes $\ell=2,3,5,\dots$ (omitting $\ell=p$) until their product exceeds $4\sqrt p$, combines the residues with the Chinese Remainder Theorem, and checks the result against Hasse's bound. The enumeration used by \texttt{ECPSEAPointCount} is limited to $p\lesssim10^6$, whereas Schoof's algorithm handles much larger primes, here $p=10^{12}+39$:

\begin{exbox}
\MatIn{83}{ECPSchoofPointCount[\{2, 3\}, 1000000000039]}

\smallskip

\MatOut{83}{1000000769568}
\end{exbox}

\noindent Here $a_p=-769528$, and the result can be checked independently by verifying that $\#E(\mathbb{F}_p)\cdot P=\mathcal{O}$ for random points $P$ with \texttt{ECPMultiplyFp}.

\medskip
\relcmds{
\cmdgroup{Point counting, group structure and point orders}{\cmd{ECPSchoofPointCount}, \cmd{ECPSEAPointCount}, \cmd{ECPGroupFp}, \cmd{ECPPointOrderFp}}
\cmdgroup{Point arithmetic and membership}{\cmd{ECPAddFp}, \cmd{ECPMultiplyFp}, \cmd{ECPPointQFp}}
\cmdgroup{Arithmetic in $\mathbb{F}_{p^r}$}{\cmd{ECPFiniteFieldInfo}, \cmd{ECPFqAdd}, \cmd{ECPFqMultiply}, \cmd{ECPFqInverse}, \cmd{ECPFqElement}, \cmd{ECPFqVector}}
\cmdgroup{Models and invariants over $\mathbb{F}_q$}{\cmd{ECPToShortWeierstrassFp}, \cmd{ECPJInvariant}}
}

\section{Curve Reductions}\label{sec:reductions}

When an elliptic curve is defined over \(\mathbb{Q}\), it is natural to ask what the curve looks like after reducing its coefficients modulo a prime \(p\). This idea is simple but one of the most useful bridges between global and local arithmetic. What actually happens is that, suppose we have an elliptic curve defined over $\mathbb{Q}$; reducing it modulo a prime $p$ produces an elliptic curve over $\mathbb{F}_{p}$. A curve may behave very nicely over \(\mathbb{Q}\), yet its reduction modulo a given prime $p$, which is an elliptic curve over $\mathbb{F}_{p}$, can reveal subtle singular behavior. In this way, reduction modulo \(p\) gives a local snapshot of the curve, and studying these snapshots prime by prime is a standard part of the arithmetic theory of elliptic curves.

\subsection{Minimal Model}
Before reducing the curve, one usually tries to choose the simplest possible model. Since the same elliptic curve can be written by many equivalent Weierstrass equations, it is important to work with a model that does not carry unnecessary powers of primes in its discriminant. This leads to the notion of a minimal Weierstrass model.

\begin{definition}[Minimal Weierstrass Model]\normalfont\cite{LozanoRobledo2011}
Let \(E/\mathbb{Q}\) be an elliptic curve. A Weierstrass equation for \(E\) is called \emph{minimal} if, among all integral Weierstrass equations for the same curve, its discriminant has minimal \(p\)-adic valuation at every prime \(p\).
\end{definition}

Once a minimal model is chosen, one reduces the coefficients modulo \(p\) and obtains a curve over \(\mathbb{F}_p\). If the reduced equation is still non-singular, then the curve has \emph{good reduction} at \(p\). If the reduced equation becomes singular, then the curve has \emph{bad reduction}. The discriminant is the key quantity here: if the discriminant does not vanish modulo \(p\), then the reduced curve remains smooth; if it does vanish, the curve acquires a singular point.

\begin{definition}[Reduction Modulo \(p\)]\normalfont\cite{Washington2008EllipticCurves}
Let
\begin{equation}
    E : y^2 = x^3 + Ax + B
\end{equation}
be an elliptic curve with \(A,B \in \mathbb{Z}\). The \emph{reduction of \(E\) modulo \(p\)} is the curve over \(\mathbb{F}_p\) obtained by reducing the coefficients \(A\) and \(B\) modulo \(p\).
\end{definition}

For a short Weierstrass equation, the discriminant is
\begin{equation}
    \Delta = -16(4A^3 + 27B^2).
\end{equation}
Thus, the first test for reduction is very direct:
\begin{equation}
    \Delta \not\equiv 0 \pmod p \quad \Longrightarrow \quad \text{good reduction at } p.
\end{equation}
 \begin{equation}
    \Delta \equiv 0 \pmod p \quad \Longrightarrow \quad \text{bad reduction at } p.
\end{equation}
In practice, for a model that is minimal at $p$, the singularity can be refined further by looking at the invariant \(c_4\).
\begin{equation}
\Delta \equiv 0 \pmod{p}
\qquad \text{but} \qquad
c_4 \not\equiv 0 \pmod{p}
\quad \Longrightarrow \quad
\text{multiplicative reduction at } p.
\end{equation}

\begin{equation}
\Delta \equiv 0 \pmod{p}
\qquad \text{and} \qquad
c_4 \equiv 0 \pmod{p}
\quad \Longrightarrow \quad
\text{additive reduction at } p.
\end{equation}
This distinction is important because it tells us not only that the curve has degenerated, but also how it has degenerated.

\begin{theorem}[Good Reduction Criterion]\normalfont\cite{Silverman2009ArithmeticEllipticCurves}
If an elliptic curve has a minimal Weierstrass model at a prime \(p\), then the reduced curve over \(\mathbb{F}_p\) is smooth if and only if the discriminant is nonzero modulo \(p\). In that case, the curve has good reduction at \(p\).
\end{theorem}

A very useful consequence of good reduction is that torsion points of order prime to \(p\) behave well under reduction. More precisely, if \(E/\mathbb{Q}\) has good reduction at \(p\) and \(m\) is not divisible by \(p\), then the reduction map gives an injective homomorphism~\cite{Washington2008EllipticCurves}
\[
E(\mathbb{Q})[m] \longrightarrow \widetilde{E}(\mathbb{F}_p).
\]
So, at good primes, finite-order rational points do not collapse unexpectedly when reduced modulo \(p\). This is one reason why reduction is so useful in explicit computations.

If the reduced curve is singular, one often refines the discussion using Kodaira symbols. In the simplest language, \(I_0\) corresponds to good reduction, \(I_n\) to multiplicative reduction, and the symbols \(II\), \(III\), \(IV\), \(I_n^*\), \(II^*\), \(III^*\), \(IV^*\) describe the various additive cases. This classification is more detailed than what is needed for a first pass, but it shows that reduction modulo \(p\) is not merely a yes-or-no question; it carries structured information about the geometry of the curve.

The reduction process can be carried out explicitly on individual curves: one first finds a minimal model, computes the discriminant, and reduces modulo a chosen prime to inspect whether the reduction is good, multiplicative, or additive.

From a computational standpoint, \texttt{ECPack} is well-suited to this topic. It can help you find a minimal model, compute the discriminant, reduce the curve modulo a chosen prime, and inspect whether the reduction is good, multiplicative, or additive. That makes it easier to move between the abstract definition and an actual curve written in code. It is also helpful when one wants to compare several primes at once, because the package can quickly show which primes are bad and which ones preserve smooth reduction. In that sense, the package turns reduction theory into something one can test and experiment with directly.

\medskip

\relcmds{
\cmdgroup{Minimal models}{\cmd{ECPIsMinimal}, \cmd{ECPIsMinimalShort}, \cmd{ECPMinimalQ}, \cmd{ECPMinimalModel}, \cmd{ECPMinimalModelShort}, \cmd{ECPFindLocalReduction}, \cmd{ECPToMinimalPoint}, \cmd{ECPFromMinimalPoint}, \cmd{ECPTransformWeierstrass}}
\cmdgroup{Discriminant and invariants}{\cmd{ECPWeierstrassDiscriminant}, \cmd{ECPWInvariants}, \cmd{ECPShortDiscriminant}, \cmd{ECPJInvariant}}
\cmdgroup{Local data: Tate's algorithm and conductor}{\cmd{ECPTate}, \cmd{ECPLocalData}, \cmd{ECPConductor}}
\cmdgroup{The reduced curve over $\mathbb{F}_p$}{\cmd{ECPToShortWeierstrassFp}, \cmd{ECPGroupFp}, \cmd{ECPSchoofPointCount}, \cmd{ECPSEAPointCount}}
}

\section{Elliptic Families and Kodaira Classification}\label{sec:kodaira}

This section, based on~\cite{Silverman1994,Tate1975,Kodaira1963II,Kodaira1963III}, provides only the introductory overview needed for the package commands. If the coefficients of a Weierstrass equation depend on a parameter \(t\), then one gets a one-parameter family of elliptic curves. The zeros of the discriminant mark the singular fibres of the family, and Kodaira's symbols provide a compact way to describe their local types. In the package, this information is used by commands that list singular fibres, report Euler numbers, and translate a fibre type into the associated gauge-algebra data when that dictionary is relevant.

\begin{table}[htbp]
\centering
\renewcommand{\arraystretch}{1.15}
\begin{tabular}{|c|c|c|c|c|}
\hline
Type & $\operatorname{ord}c_4$ & $\operatorname{ord}c_6$ & $\operatorname{ord}\Delta$ & Euler number\\
\hline
$I_0$ & $\ge0$ & $\ge0$ & $0$ & $0$\\
$I_n$ $(n\ge1)$ & $0$ & $0$ & $n$ & $n$\\
$II$ & $\ge1$ & $1$ & $2$ & $2$\\
$III$ & $1$ & $\ge2$ & $3$ & $3$\\
$IV$ & $\ge2$ & $2$ & $4$ & $4$\\
$I_0^{*}$ & $\ge2$ & $\ge3$ & $6$ & $6$\\
$I_n^{*}$ $(n\ge1)$ & $2$ & $3$ & $n+6$ & $n+6$\\
$IV^{*}$ & $\ge3$ & $4$ & $8$ & $8$\\
$III^{*}$ & $3$ & $\ge5$ & $9$ & $9$\\
$II^{*}$ & $\ge4$ & $5$ & $10$ & $10$\\
\hline
\end{tabular}
\caption{Kodaira's table: the type of a singular fibre from the orders of vanishing of $c_4$, $c_6$ and $\Delta$ at a point where the Weierstrass equation is minimal.}
\label{tab:kodaira}
\end{table}

A minimal model at a point \(t_0\) is detected by the vanishing orders of the invariants \(c_4\), \(c_6\), and \(\Delta\). For a first reading, it is enough to remember that the fibre is smooth when \(\Delta(t_0)\neq 0\), and that the symbols \(I_n\), \(II\), \(III\), \(IV\), \(I_n^*\), \(II^*\), \(III^*\), and \(IV^*\) are the standard names used by the commands in this package. The same local data can also be viewed at \(t=\infty\) by passing to \(s=1/t\), which is why the package can treat families over \(\mathbb{P}^1\) uniformly.

The related command \texttt{ECPFibreList} returns the singular fibres and their Euler numbers, while \texttt{ECPFibreAt} inspects a fibre at a chosen parameter value and \texttt{ECPGaugeAlgebra} records the corresponding gauge-algebra dictionary when needed.

\relcmds{
\cmdgroup{Kodaira classification}{\cmd{ECPKodairaType}, \cmd{ECPFibreAt}, \cmd{ECPFibreAtInfinity}, \cmd{ECPFibreList}}
\cmdgroup{Miscellaneous}{\cmd{ECPJFunction}, \cmd{ECPCoulombDimension}, \cmd{ECPGaugeAlgebra}, \cmd{ECPTateForm}, \cmd{ECPTorsionFamily}}
}

\newpage
\section{Cryptography and Algebraic Pairings}\label{sec:crypto}

The two preceding sections studied the finite groups $E(\mathbb{F}_q)$ and the curves obtained from a curve over $\mathbb{Q}$ by reduction modulo $p$. These finite groups are also the setting of elliptic-curve cryptography, which uses the group law over $\mathbb{F}_q$ directly, and of pairings, which relate the additive group $E(\mathbb{F}_q)$ to the multiplicative group $\mathbb{F}_q^{\times}$. In this section we first recall the discrete logarithm problem on an elliptic curve and the ElGamal protocol, which \texttt{ECPack} implements over prime fields. We then describe Miller's algorithm and the Weil and Tate pairings, which are implemented over $\mathbb{F}_p$ and $\mathbb{F}_{p^r}$, together with a routine that produces a basis of the full $n$-torsion subgroup.

\subsection{Elliptic Curve Cryptography and ElGamal Encryption}\label{subsec:elgamal}

Elliptic curves became especially important in modern cryptography because they offer a very efficient group law, while the inverse problem remains difficult in practice. Over a finite field, point addition and scalar multiplication can be performed quickly, but recovering the scalar from a given multiple of a point is believed to be computationally hard. This asymmetry is the starting point for elliptic curve cryptography (ECC)~\cite{Koblitz1987,Miller1985,Washington2008EllipticCurves,HMV2004}.

A central problem here is the \emph{Elliptic Curve Discrete Logarithm Problem} (ECDLP). Suppose \(P\) is a point on an elliptic curve and \(Q = nP\) for some unknown integer \(n\). The ECDLP asks for the recovery of \(n\) from \(P\) and \(Q\). For suitably chosen curves, no efficient general method is known for solving this problem~\cite{Washington2008EllipticCurves}, which is why elliptic curves are useful for secure communication. In contrast, the forward operation \(P \mapsto nP\) is easy to compute, especially with double-and-add methods and related fast scalar-multiplication algorithms.

A standard application is the elliptic-curve version of the ElGamal protocol~\cite{ElGamal1985,Koblitz1987}. In this setting, a public point \(P\) is chosen on the curve, together with a public point \(Q = dP\), where \(d\) is the secret key. To encrypt a message, one first encodes the message as a point \(M\) on the curve and then chooses a random integer \(k\). The ciphertext is the pair
\begin{equation}
    C_1 = kP, \qquad C_2 = M + kQ.
\end{equation}
Since \(Q = dP\), the legitimate receiver can recover the message by computing
\begin{equation}
    M = C_2 - dC_1.
\end{equation}
The security of this procedure depends on the difficulty of recovering \(d\) from \(P\) and \(Q\). The group law makes encryption and decryption efficient, while the underlying discrete logarithm problem makes unauthorized recovery difficult. \texttt{ECPack} implements Koblitz's embedding of an integer message as a point on the curve~\cite{Koblitz1987} together with ElGamal encryption and decryption over prime fields. The parameters used in the manual are toy values: the routines are meant for experimentation with the arithmetic, not for secure communication.

\subsection{Pairings and Miller's Algorithm}\label{subsec:pairings}

In more advanced cryptographic constructions, one needs a bilinear pairing. A pairing takes two points of an elliptic curve and returns an element of the multiplicative group of the field, so that the additive group law of the curve is converted into a multiplicative one in a controlled way. Throughout this subsection $E$ is defined over $\mathbb{F}_q$, $q=p^r$, and $n$ is a positive integer prime to $p$. Over an algebraic closure, $E[n]\cong(\mathbb{Z}/n\mathbb{Z})^2$, and we write $\mu_n$ for the group of $n$-th roots of unity.

\paragraph{The Weil Pairing.} The most classical pairing is the Weil pairing.

\begin{definition}[Weil Pairing]\normalfont\cite{Washington2008EllipticCurves}
The Weil pairing is a map
\begin{equation}
    e_n : E[n] \times E[n] \to \mu_n
\end{equation}
with the following properties:
\begin{enumerate}
    \item[(a)] \emph{bilinearity:} $e_n(aP,bQ) = e_n(P,Q)^{ab}$ for integers $a,b$ and $P,Q\in E[n]$;
    \item[(b)] \emph{alternating:} $e_n(P,P)=1$, and hence $e_n(Q,P)=e_n(P,Q)^{-1}$;
    \item[(c)] \emph{non-degeneracy:} if $e_n(P,Q)=1$ for all $Q\in E[n]$, then $P=\mathcal{O}$.
\end{enumerate}
\end{definition}

If $P,Q$ is a basis of $E[n]$, then $e_n(P,Q)$ is a primitive $n$-th root of unity. In particular, if the full $n$-torsion subgroup is contained in $E(\mathbb{F}_q)$, then $\mu_n\subset\mathbb{F}_q$, that is, $n$ divides $q-1$.

\paragraph{Miller's Function and Miller's Algorithm\normalfont~\cite{Miller2004}.} The Weil pairing is computed with Miller's algorithm, whose basic object is a rational function on $E$. For a point $P$ and an integer $n\ge 1$, let $f_{n,P}$ be a rational function with divisor
\begin{equation}
    \operatorname{div}\,f_{n,P}=n(P)-([n]P)-(n-1)(\mathcal{O}),
\end{equation}
normalised as in Miller's algorithm; for $P\in E[n]$ this divisor is $n(P)-n(\mathcal{O})$. Let $\ell_{T,R}$ denote the line through $T$ and $R$ (the tangent if $T=R$) and $v_S$ the vertical line through $S$. Comparing divisors, one finds
\begin{equation}
    f_{i+j,P}=f_{i,P}\,f_{j,P}\,\frac{\ell_{[i]P,\,[j]P}}{v_{[i+j]P}},
    \label{eq:miller_recursion}
\end{equation}
which is the geometric chord-and-tangent construction of Section~\ref{sec:rationals} in the language of divisors. Miller's algorithm evaluates $f_{n,P}(Q)$ at a point $Q\neq\mathcal{O}$ by applying \eqref{eq:miller_recursion} along the binary expansion of $n$ (double-and-add). At each step, only the values of the line and of the vertical line at $Q$ are computed, so the function itself is never written down, and the cost is $O(\log n)$ operations in $\mathbb{F}_q$~\cite{Miller2004}. The evaluation fails if $Q$ is a zero or a pole of one of the intermediate functions, for instance when $Q$ lies in the subgroup generated by $P$.

For $P,Q\in E[n]$ the Weil pairing is given by
\begin{equation}
    e_n(P,Q)=(-1)^n\,\frac{f_{n,P}(Q)}{f_{n,Q}(P)} .
    \label{eq:weil_miller}
\end{equation}
Both values in \eqref{eq:weil_miller} are defined and nonzero exactly when $P$ and $Q$ are independent; if one of them lies in the cyclic group generated by the other, the pairing equals $1$.

\paragraph{The Tate Pairing\normalfont~\cite{Lichtenbaum1969,FreyRueck1994,Washington2008EllipticCurves}.} The second pairing in common use is the Tate pairing. It is defined on $E(\mathbb{F}_q)[n]\times E(\mathbb{F}_q)/nE(\mathbb{F}_q)$ and takes values in $\mathbb{F}_q^{\times}/(\mathbb{F}_q^{\times})^n$. If $n$ divides $q-1$, one makes the value unique by raising it to the power $(q-1)/n$, which gives the \emph{reduced Tate pairing}
\begin{equation}
    t_n(P,Q)=f_{n,P}(Q)^{(q-1)/n}\in\mu_n .
\end{equation}
It is bilinear and non-degenerate but, unlike the Weil pairing, it is not alternating, so $t_n(P,P)$ need not equal $1$. Its evaluation requires only one Miller function rather than two. The condition $n\mid q-1$ is the reason why one often has to pass to an extension field: the smallest $k$ with $n\mid q^k-1$ is called the embedding degree of $n$~\cite{MOV1993,FreyRueck1994}.

\paragraph{Implementation in \texttt{ECPack}.} \texttt{ECPack} evaluates the Miller function with \texttt{ECPMillerFp}, and computes the Weil pairing \eqref{eq:weil_miller} and the reduced Tate pairing with \texttt{ECPWeilPairingFp} and \texttt{ECPTatePairingFp}. All three work over $\mathbb{F}_p$ and over $\mathbb{F}_{p^r}$, in every characteristic for which the pairing is defined, and for a curve in any of the accepted forms. The routine \texttt{ECPTorsionBasisFp} finds a basis $\{P,Q\}$ of $E[n]$ inside $E(\mathbb{F}_q)$, with $P$ of exact order $n$ and $e_n(P,Q)$ a primitive $n$-th root of unity; this requires $E[n]\subset E(\mathbb{F}_q)$. For the curve $y^2=x^3+2$ over $\mathbb{F}_7$, one has $9$ points and $E(\mathbb{F}_7)\cong\mathbb{Z}/3\times\mathbb{Z}/3$, so $E[3]\subset E(\mathbb{F}_7)$ and $3\mid 7-1$:

\begin{exbox}
\MatIn{80}{ECPTorsionBasisFp[3, \{0, 2\}, 7]}

\smallskip

\MatOut{80}{\{\{0, 4\}, \{5, 1\}\}}

\medskip

\MatIn{81}{ECPWeilPairingFp[3, \{0, 4\}, \{5, 1\}, \{0, 2\}, 7]}

\smallskip

\MatOut{81}{2}

\medskip

\MatIn{82}{\{ECPWeilPairingFp[3, \{0, 4\}, \{0, 4\}, \{0, 2\}, 7], ECPWeilPairingFp[3, \{5, 1\}, \{0, 4\}, \{0, 2\}, 7]\}}

\smallskip

\MatOut{82}{\{1, 4\}}
\end{exbox}

\noindent Since $2^3\equiv1 \pmod 7$ and $2\not\equiv1$, the value $e_3(P,Q)=2$ is a primitive cube root of unity, as it must be for a basis. The pairing is alternating, $e_3(P,P)=1$, and $e_3(Q,P)=4=2^{-1}$ modulo $7$.

\relcmds{
\cmdgroup{Group law over $\mathbb{F}_p$}{\cmd{ECPAddFp}, \cmd{ECPMultiplyFp}, \cmd{ECPPointOrderFp}, \cmd{ECPGroupFp}}
\cmdgroup{Message embedding and ElGamal}{\cmd{ECPFindPointFp}, \cmd{ECPElGamalEncodeFp}, \cmd{ECPElGamalDecodeFp}}
\cmdgroup{Miller's algorithm and pairings}{\cmd{ECPMillerFp}, \cmd{ECPWeilPairingFp}, \cmd{ECPTatePairingFp}, \cmd{ECPTorsionBasisFp}}
}

\section{Future Directions}\label{sec:future}

We plan to extend \texttt{ECPack.wl} in the following directions.

\begin{enumerate}
\item \textbf{Torsion via Doud's Algorithm.} \texttt{ECPFindTorsionPoints} uses the Nagell--Lutz theorem, which requires factoring the discriminant and falls back to a bounded search for large $|\Delta|$. Doud's algorithm~\cite{Doud1998} instead bounds the torsion order by point counts modulo primes of good reduction and tests the $n$-division points of the analytic parametrization by $\wp$ for rationality; we intend to implement it, with exact verification of the candidates.
\item \textbf{Sharper Rank Bounds.} The bounds of Section~\ref{subsec:primebound} need a rational point of order $2$. Computing the $2$-Selmer group gives $R_E\le\dim_{\mathbb{F}_2}\mathrm{Sel}^{(2)}(E/\mathbb{Q})-\dim_{\mathbb{F}_2}E(\mathbb{Q})[2]$ for every curve over $\mathbb{Q}$, with equality exactly when $\mathrm{Sha}(E/\mathbb{Q})[2]=0$~\cite{Selmer1951DiophantineEquation,Cremona1997AlgorithmsModularEllipticCurves}.
\item \textbf{Elliptic Curves over $\mathbb{C}$.} We plan commands for the lattice $\Lambda$, the periods, $g_2$, $g_3$, the modular parameter $\tau$ and the uniformization $\mathbb{C}/\Lambda\to E(\mathbb{C})$, which is the case relevant to elliptic Feynman integrals and is also needed in item~1.
\item \textbf{Further Extensions.} A genuine Schoof--Elkies--Atkin point count (\texttt{ECPSEAPointCount} is currently an enumeration), canonical heights through local heights~\cite{Silverman1988ComputingHeights}, the Tamagawa numbers and real period for a numerical check of the BSD formula, isogenies (V\'elu), the Shioda--Tate formula for families, and further cryptographic protocols such as ECDH and ECDSA.
\end{enumerate}

\section*{Acknowledgements}
During the preparation of this work, we have used the AI assistants Claude (Anthropic) and ChatGPT (OpenAI) to help refine and test the Wolfram Language code of the package and to check test cases. All mathematical content, the package design, and the final code, examples, and text were reviewed by the authors.


\appendix
\sloppy

\section{Characteristic and Finite Fields}\label{app:fields}

The main text refers to the characteristic of a field and to the fields $\mathbb{F}_{p^r}$ without defining them, although these notions decide which of the transformations and algorithms of \texttt{ECPack} apply to a given curve. This appendix collects the definitions and explains how $\mathbb{F}_{p^r}$ is represented inside the package.

\subsection{Characteristic}\label{app:char}
The characteristic of a field $K$ is the smallest $n\ge 1$ with $n\cdot 1=0$, or $0$ if no such $n$
exists. It is always $0$ or a prime: if $n=ab$ with $a,b>1$, then $(a\cdot1)(b\cdot1)=0$ would give
zero divisors. Thus $\operatorname{char}\mathbb{Q}=\operatorname{char}\mathbb{R}=\operatorname{char}\mathbb{C}=0$
and $\operatorname{char}\mathbb{F}_{p^r}=p$. The short Weierstrass form
$y^2=x^3+Ax+B$ requires $\operatorname{char}K\neq2,3$ (completing the square and the cube divides by $2$
and $3$); in characteristic $2$ and $3$ one must use the normal forms of \texttt{ECPNormalFormCharTwoThree}
and \texttt{ECPToShortWeierstrassFp}.

\subsection{Prime Fields and Prime Field Extensions}\label{app:prime}
For a prime $p$, $\mathbb{F}_p=\mathbb{Z}/p\mathbb{Z}$ is a field, while $\mathbb{Z}/n\mathbb{Z}$ is not for
composite $n$ (e.g.\ $2\cdot3\equiv0 \bmod 6$). Over $\mathbb{F}_p$ all coefficients, coordinates and group-law
computations are carried out modulo $p$.

For $q=p^r$ with $r\ge 2$, $\mathbb{F}_q$ is \emph{not} $\mathbb{Z}/q\mathbb{Z}$. It is constructed as
$\mathbb{F}_p[x]/(f)$ for any irreducible $f\in\mathbb{F}_p[x]$ of degree $r$, and it is unique up to
isomorphism; $\mathbb{F}_q\subseteq\mathbb{F}_{q^n}$ for all $n\ge1$, and $\mathbb{F}_q^\times$ is cyclic.
For example, $\omega^2+\omega+1$ has no root in $\mathbb{F}_2$, so
\begin{equation}
  \mathbb{F}_4=\mathbb{F}_2[\omega]/(\omega^2+\omega+1)=\{0,1,\omega,\omega^2\},\qquad \omega^2=\omega+1 .
\end{equation}
ECPack represents $\mathbb{F}_{p^r}$ as $\mathbb{F}_p[x]/(f)$ and writes nonzero elements as powers $g^k$ of a primitive
element $g$, so that multiplication is addition of exponents modulo $q-1$ (see \texttt{ECPFiniteFieldInfo}).
Since $E(\mathbb{F}_q)\subseteq E(\mathbb{F}_{q^n})$ and both are finite, every curve over $\mathbb{F}_q$ has a finite
group of points over each $\mathbb{F}_{q^n}$.

\paragraph*{Relevant ECPack commands.}\ 
\texttt{ECPNormalFormCharTwoThree}, \texttt{ECPCubicToWeierstrassCharTwoThree}, \texttt{ECPToShortWeierstrassFp},
\texttt{ECPJInvariant}, \texttt{ECPFiniteFieldInfo}, \texttt{ECPFqAdd}, \texttt{ECPFqMultiply},
\texttt{ECPFqInverse}, \texttt{ECPFqElement}, \texttt{ECPFqVector}, \texttt{ECPGroupFp}, \texttt{ECPAddFp},
\texttt{ECPPointQFp}.

\section{The Frobenius Endomorphism}\label{app:frobenius}

Let $E/\mathbb{F}_q$ be given by $y^2=x^3+Ax+B$ with $A,B\in\mathbb{F}_q$, $q=p^r$, and let $\overline{\mathbb{F}}_q$
be an algebraic closure. An \emph{endomorphism} of $E$ is a morphism $E\to E$ fixing $\mathcal{O}$ (hence a group
homomorphism). For example, the multiplication maps $[n]\colon P\mapsto nP$ are endomorphisms, with $\ker[n]=E[n]$.

\subsection{Frobenius on the Field and on the Curve}\label{app:frob-def}
In characteristic $p$ the map $a\mapsto a^p$ is a ring homomorphism, since $(a+b)^p=a^p+b^p$ (the binomial
coefficients $\binom{p}{k}$, $0<k<p$, vanish mod $p$). Its iterate $a\mapsto a^q$ fixes exactly $\mathbb{F}_q$:
\begin{equation}
  \mathbb{F}_q=\{a\in\overline{\mathbb{F}}_q : a^q=a\}.
\end{equation}
The $q$-power Frobenius of the curve is
\begin{equation}
  \phi_q\colon E\to E,\qquad (x,y)\mapsto(x^q,y^q),\qquad \phi_q(\mathcal{O})=\mathcal{O}.
\end{equation}
It maps $E$ to itself: raising $y^2=x^3+Ax+B$ to the $q$-th power and using $A^q=A$, $B^q=B$ gives
$(y^q)^2=(x^q)^3+Ax^q+B$. The group law is given by rational functions with coefficients in $\mathbb{F}_q$, which
commute with $\phi_q$, so $\phi_q(P+Q)=\phi_q(P)+\phi_q(Q)$ and $\phi_q$ is an endomorphism. By the characterization
above,
\begin{equation}
  E(\mathbb{F}_{q^n})=\{P\in E(\overline{\mathbb{F}}_q):\ \phi_q^{\,n}(P)=P\}=\ker(\phi_q^{\,n}-1).
  \label{eq:frobfixed}
\end{equation}
The Frobenius endomorphism on finite-field points can be used to describe fixed points over extension fields and to relate point counts over \(\mathbb{F}_{q^n}\) to the kernel of \(\phi_q^n-1\).

\subsection{Trace of Frobenius and Point Counts}\label{app:trace}
\begin{theorem}\label{thm:trace}
Let $a_q=q+1-\#E(\mathbb{F}_q)$.
\begin{enumerate}
  \item[(a)] $\phi_q^2-a_q\phi_q+q=0$ in $\operatorname{End}(E)$.
  \item[(b)] Writing $X^2-a_qX+q=(X-\alpha)(X-\beta)$, one has $|\alpha|=|\beta|=\sqrt q$ and
  $\#E(\mathbb{F}_{q^n})=q^n+1-\alpha^n-\beta^n$ for all $n\ge1$. In particular $|a_q|\le2\sqrt q$ (Hasse).
\end{enumerate}
\end{theorem}
\noindent\emph{Sketch.} $\phi_q$ has degree $q$ and $\phi_q-1$ is separable, so by
\eqref{eq:frobfixed}, $\#E(\mathbb{F}_q)=\deg(\phi_q-1)=(\phi_q-1)(\hat\phi_q-1)=1-(\phi_q+\hat\phi_q)+q$, with
$\hat\phi_q$ the dual isogeny. Hence $a_q=\phi_q+\hat\phi_q$ and $\phi_q\hat\phi_q=q$, which gives (a);
(b) follows from $\#E(\mathbb{F}_{q^n})=\deg(\phi_q^n-1)$ and the Weil bound on $\alpha$. See
Silverman~\cite{Silverman2009ArithmeticEllipticCurves} (Ch.~V) or Washington~\cite{Washington2008EllipticCurves} (Ch.~4).

\subsection{Action on Torsion and the Tate Module}\label{app:tate}
For $\ell\ne p$ prime, $E[\ell^n]\cong(\mathbb{Z}/\ell^n\mathbb{Z})^2$, and $\phi_q$ acts on each $E[\ell^n]$
compatibly. The Tate module
\begin{equation}
  T_\ell(E)=\varprojlim_n E[\ell^n]\cong\mathbb{Z}_\ell^{\,2},\qquad V_\ell(E)=T_\ell(E)\otimes_{\mathbb{Z}_\ell}\mathbb{Q}_\ell\cong\mathbb{Q}_\ell^{\,2},
\end{equation}
organises these actions into a single $2\times2$ matrix over $\mathbb{Q}_\ell$ whose characteristic polynomial is
$X^2-a_qX+q$, i.e.\ $\operatorname{tr}\phi_q=a_q$, $\det\phi_q=q$. Restricting to $E[\ell]$ gives
$a_q\bmod\ell$, which is the quantity computed by Schoof's algorithm in Section~\ref{sec:schoof}; $a_q$ is then recovered by the Chinese Remainder Theorem together with Hasse's bound.

\paragraph*{Relevant ECPack commands.}\ 
\texttt{ECPGroupFp}, \texttt{ECPAddFp}, \texttt{ECPMultiplyFp}, \texttt{ECPPointOrderFp},
\texttt{ECPSchoofPointCount}, \texttt{ECPSEAPointCount}.

\section{User Manual}\label{app:manual}
All public commands of the package \texttt{ECPack} begin with the prefix \texttt{ECP} (variables with \texttt{\$ECP}), so that all of them can be listed with \texttt{Names["ECPack`ECP*"]}. The package is loaded with \texttt{Get["ECPack.wl"]}.
Each entry below provides the call forms, arguments, returned value, remarks, and one example in the style of a Mathematica notebook. Unless stated otherwise, \texttt{curve} may be a long model $\{a_1,a_2,a_3,a_4,a_6\}$, a short pair $\{A,B\}$, an equation, a polynomial in $x,y$, a curve Association such as \texttt{\$ECPCurrentEC}, or \texttt{Automatic} (the session curve). The examples use the curves 37a1: $y^2+y=x^3-x$ (rank 1, conductor 37), $y^2=x^3-36x$ (CM, torsion $\mathbb Z/2\times\mathbb Z/2$) and $y^2=x^3+73$.
\subsection{Index of Commands}
Commands marked \slowmark\ can take a long time to compute for some inputs; the reason is given in their entries below.
\begin{description}\setlength{\itemsep}{0pt}

\item[\texttt{ECPAdd}] Add two rational points on a curve.
\item[\texttt{ECPAddFp}] Add two points over $\mathbb F_p$ or $\mathbb F_{p^r}$ (short form: odd $p$; long form: any $p$).
\item[\texttt{ECPAnalyticRank}\slowmark] Analytic rank, root number and conductor, assuming BSD.
\item[\texttt{ECPAreProjectivelyEquivalent}] Test whether two projective points are proportional.
\item[\texttt{ECPBachetDuplication}] Bachet's duplication formula: $x(2P)$ as a rational function of $x(P)$.
\item[\texttt{ECPCanonicalHeight}\slowmark] Canonical (N\'eron--Tate) height of a point.
\item[\texttt{ECPChordAndTangent}] Third intersection point of the tangent line at $P$ with the cubic.
\item[\texttt{ECPComposeTransformations}] Compose two Weierstrass isomorphisms $\{u,r,s,t\}$.
\item[\texttt{ECPConductor}\slowmark] Conductor $N$ of an integral Weierstrass model.
\item[\texttt{ECPCoulombDimension}] Scaling dimension $12/(12-\operatorname{ord}\Delta)$ of a Kodaira singularity.
\item[\texttt{ECPCubicToWeierstrass}] Reduce a plane cubic with a known rational point to Weierstrass form (characteristic 0 or $>3$).
\item[\texttt{ECPCubicToWeierstrassCharTwoThree}] Same, given an inflection point, in characteristic 2 or 3 (also over $\mathbb F_{p^r}$).
\item[\texttt{ECPElGamalDecodeFp}] Decrypt an ElGamal ciphertext over $\mathbb F_p$ (prime fields only).
\item[\texttt{ECPElGamalEncodeFp}] Encrypt a message by elliptic-curve ElGamal over $\mathbb F_p$ (prime fields only).
\item[\texttt{ECPExtractCoefficients}] Read $\{a_1,a_2,a_3,a_4,a_6\}$ off a curve in any form.
\item[\texttt{ECPFibreAt}] Kodaira fibre of a family over $\mathbb Q(t)$ at a point $t_0$.
\item[\texttt{ECPFibreAtInfinity}] Kodaira fibre at $t=\infty$ of a family over $\mathbb P^1$.
\item[\texttt{ECPFibreList}] All singular fibres of a family over $\mathbb P^1$, with the Euler-number check.
\item[\texttt{ECPFindGenerators}\slowmark] Heuristic search for independent generators of the Mordell--Weil group.
\item[\texttt{ECPFindLocalReduction}] Local Laska--Kraus reduction step at a prime $p$.
\item[\texttt{ECPFindPointFp}] Koblitz-embed a message as a curve point over $\mathbb F_p$ (prime fields only).
\item[\texttt{ECPFindPointOrder}] Order of a point by bounded search (Mazur's bound by default).
\item[\texttt{ECPFindPointOrderFp}\slowmark] Same as \texttt{ECPPointOrderFp}: exact order of a point over $\mathbb F_p$ or $\mathbb F_{p^r}$.
\item[\texttt{ECPFindTorsionPoints}\slowmark] All rational torsion points (Nagell--Lutz).
\item[\texttt{ECPFiniteFieldInfo}] Describe the field polynomial and generator used for $\mathbb F_{p^r}$.
\item[\texttt{ECPFormattedEquation}] Display the long Weierstrass equation.
\item[\texttt{ECPFqAdd}] Add two elements of $\mathbb F_{p^r}$.
\item[\texttt{ECPFqElement}] Build a field element from its coordinate vector.
\item[\texttt{ECPFqInverse}] Multiplicative inverse in $\mathbb F_{p^r}$.
\item[\texttt{ECPFqMultiply}] Multiply two elements of $\mathbb F_{p^r}$.
\item[\texttt{ECPFqVector}] Coordinate vector of a field element.
\item[\texttt{ECPFromMinimalPoint}] Inverse of \texttt{ECPToMinimalPoint}.
\item[\texttt{ECPFromShortCoordinates}] Map a point on the short model back to the long model.
\item[\texttt{ECPFromShortCoordinatesFp}] Finite-field analogue of \texttt{ECPFromShortCoordinates}.
\item[\texttt{ECPGaugeAlgebra}] F-theory gauge algebra of a Kodaira type (split or non-split).
\item[\texttt{ECPGroupFp}\slowmark] The group $E(\mathbb F_q)$ of points over $\mathbb F_p$ or $\mathbb F_{p^r}$.
\item[\texttt{ECPHeightMatrix}\slowmark] Pairwise height-pairing matrix (its determinant is the regulator).
\item[\texttt{ECPHeightPairing}\slowmark] N\'eron--Tate height pairing of two points.
\item[\texttt{ECPIndependentQ}\slowmark] Test independence of points via the height-matrix determinant.
\item[\texttt{ECPIsMinimal}] Test true Weierstrass minimality of a long model.
\item[\texttt{ECPIsMinimalShort}] Test minimality of a short model (or dispatch to \texttt{ECPIsMinimal}).
\item[\texttt{ECPJFunction}] $j$-invariant $c_4^3/\Delta$ of a family as a rational function.
\item[\texttt{ECPJInvariant}] $j$-invariant over $\mathbb Q$, $\mathbb F_p$ or $\mathbb F_{p^r}$.
\item[\texttt{ECPKodairaType}] Kodaira type and Euler number from the orders of $c_4,c_6,\Delta$.
\item[\texttt{ECPLegendreToWeierstrass}] Convert Legendre form to short Weierstrass form.
\item[\texttt{ECPLocalData}\slowmark] Tate's-algorithm data at every bad prime.
\item[\texttt{ECPMillerFp}] Miller's function $f_{n,P}$ evaluated at a point $Q$ over $\mathbb F_p$ or $\mathbb F_{p^r}$.
\item[\texttt{ECPMinimalModel}\slowmark] Global minimal (long) Weierstrass model.
\item[\texttt{ECPMinimalModelShort}] Minimal model for a short Weierstrass pair.
\item[\texttt{ECPMinimalQ}] Test minimality at one prime or globally.
\item[\texttt{ECPMultiply}] Scalar multiplication of a rational point.
\item[\texttt{ECPMultiplyFp}] Scalar multiplication over $\mathbb F_p$ or $\mathbb F_{p^r}$.
\item[\texttt{ECPNaiveHeight}] Naive logarithmic height of a point.
\item[\texttt{ECPNormalFormCharTwoThree}] Normal form of a long model in characteristic 2 or 3 (also over $\mathbb F_{p^r}$).
\item[\texttt{ECPOnCurveQ}] Test whether a point lies on a given curve.
\item[\texttt{ECPParseAffineCurve}] Normalise a curve to $\{\text{polynomial},\{x,y\}\}$.
\item[\texttt{ECPParseHomogeneousCurve}] Homogenise a cubic with respect to $z$.
\item[\texttt{ECPPointOrder}] Order of a rational point.
\item[\texttt{ECPPointOrderFp}\slowmark] Order of a point over $\mathbb F_p$ or $\mathbb F_{p^r}$ via Lagrange's theorem.
\item[\texttt{ECPPointQ}] Membership test for a point on a curve (coordinates as two arguments).
\item[\texttt{ECPPointQFp}] Membership test over $\mathbb F_p$ or $\mathbb F_{p^r}$.
\item[\texttt{ECPPolynomialHomogeneous}] Exact homogeneity test.
\item[\texttt{ECPPolynomialTransform}] Linear change of variables in a polynomial.
\item[\texttt{ECPPrimeDivisorCount}] Number $\nu(n)$ of distinct primes dividing $n$, or dividing $\Delta$ of a curve.
\item[\texttt{ECPQuarticToWeierstrass}] Convert $y^2=$ quartic to short Weierstrass form.
\item[\texttt{ECPRankBound}\slowmark] Proven upper bound for the rank of a curve with a rational point of order 2.
\item[\texttt{ECPRankUpperBound}\slowmark] Heuristic (non-rigorous) rank estimate from the discriminant.
\item[\texttt{ECPSchoofPointCount}\slowmark] Number of points over $\mathbb F_p$ by Schoof's algorithm (large $p$).
\item[\texttt{ECPSEAPointCount}\slowmark] Exact point count by enumeration over $\mathbb F_p$ ($p\lesssim10^6$) or $\mathbb F_{p^r}$.
\item[\texttt{ECPSearchPoints}] Bounded rational-point search.
\item[\texttt{ECPSetEC}] Set the global session curve.
\item[\texttt{ECPShortDiscriminant}] Discriminant $-16(4A^3+27B^2)$ of a short model.
\item[\texttt{ECPShortWeierstrassForm}] Gives short form curve with integer coefficients.

\item[\texttt{ECPTate}] Tate's algorithm at one prime.
\item[\texttt{ECPTateForm}] Coefficients with prescribed vanishing orders at $t=0$ (F-theory Tate form).
\item[\texttt{ECPTatePairingFp}\slowmark] Reduced Tate pairing $t_n(P,Q)$ over $\mathbb F_p$ or $\mathbb F_{p^r}$.
\item[\texttt{ECPToLongCoordinates}] Map coordinates of a transformed model back to the long model.
\item[\texttt{ECPToMinimalPoint}] Map a point under a Weierstrass isomorphism $\{u,r,s,t\}$.
\item[\texttt{ECPTorsionBasisFp}\slowmark] Basis of the full $n$-torsion subgroup $E[n]\subset E(\mathbb F_q)$.
\item[\texttt{ECPTorsionFamily}] Universal curve with a rational $N$-torsion point (Tate--Kubert form), $N=4,\dots,7$.
\item[\texttt{ECPTorsionSubgroup}\slowmark] Torsion subgroup (order, structure, generators, points).
\item[\texttt{ECPToShortCoordinates}] Map a point on the long model to the short model.
\item[\texttt{ECPToShortCoordinatesFp}] Finite-field analogue of \texttt{ECPToShortCoordinates} (prime field).
\item[\texttt{ECPToShortWeierstrass}] Reduce a long model to short Weierstrass form (characteristic 0 or $>3$).
\item[\texttt{ECPToShortWeierstrassFp}] Reduce a long model to its simplest form over $\mathbb F_p$ or $\mathbb F_{p^r}$ (any $p$).
\item[\texttt{ECPToWeierstrassForm}] Read $\{a_1,\dots,a_6\}$ off a cubic already in Weierstrass shape.
\item[\texttt{ECPTransformWeierstrass}] Apply the isomorphism $\{u,r,s,t\}$ to the coefficients.
\item[\texttt{ECPWeierstrassDiscriminant}] Full discriminant $\Delta$ of a long model.
\item[\texttt{ECPWeilPairingFp}] Weil pairing $e_n(P,Q)$ over $\mathbb F_p$ or $\mathbb F_{p^r}$.
\item[\texttt{ECPWInvariants}] The invariants $b_2,b_4,b_6,b_8,c_4,c_6$ and $\Delta$.
\end{description}
\noindent\texttt{ECPShortWierstrass} is kept as a legacy alias of \texttt{ECPShortWeierstrassForm}. The package also defines the variables \texttt{\$ECPVersion} (the version string, \texttt{"1.0"}), \texttt{\$ECPCurrentEC} (the session curve set by \texttt{ECPSetEC}), \texttt{\$ECPNagellLutzDiscLimit} (the largest $|\Delta|$ for exact Nagell--Lutz in \texttt{ECPFindTorsionPoints}) and \texttt{\$ECPFqTableLimit} (the largest $q$ for which the arithmetic tables of $\mathbb F_q$ are built).

\subsection{Function Descriptions}

\fname{ECPAdd}
\usageline{ECPAdd[P1, P2, curve]}
\args\ \texttt{P1}, \texttt{P2}: points $\{x,y\}$ or \texttt{Infinity} in the coordinates of \texttt{curve}; \texttt{curve}: any curve form.
\returns\ The point $P_1+P_2$ in the same coordinates as \texttt{curve}, as $\{x,y\}$ or \texttt{Infinity}.
\notes\ Both points and the curve are validated once at the entry point; the core uses the doubling formula if $P_1=P_2$ and the chord formula otherwise (a vertical chord gives \texttt{Infinity}).

\begin{exbox}
\textbf{Example.} (0,0) + (1,0) on 37a1

\MatIn{4}{ECPAdd[\{0, 0\}, \{1, 0\}, \{0, 0, 1, -1, 0\}]}

\smallskip

\MatOut{4}{\{-1, -1\}}
\end{exbox}

\fname{ECPAddFp}
\usageline{ECPAddFp[P1, P2, a, p]\quad ECPAddFp[P1, P2, curve, p]\quad ECPAddFp[P1, P2, curve, \{p, r\}]}
\args\ \texttt{P1}, \texttt{P2}: points or \texttt{Infinity} (integer coordinates need not be reduced mod $p$). \texttt{a}: the $x$-coefficient of $y^2=x^3+ax+b$ (short overload, odd $p$ only; $b$ is not needed, so membership is not tested). \texttt{curve}: any curve form (long overload; it checks that both points lie on the curve; every prime allowed). \texttt{\{p,r\}}: extension field, $r\ge 2$.
\returns\ The sum over $\mathbb F_p$ with coordinates in $[0,p)$, or \texttt{Infinity}. Over $\mathbb F_{p^r}$ coordinates are written $0$ or $g^k$.
\notes\ Points that do not lie on one common curve mod $p$ raise a message and return \texttt{\$Failed}. Over an extension field give the full curve and field elements in the notation of \texttt{ECPFiniteFieldInfo}.

\begin{exbox}
\textbf{Example.} (5,1) + (6,3)

\MatIn{5}{ECPAddFp[\{5, 1\}, \{6, 3\}, 2, 17]}

\smallskip

\MatOut{5}{\{10, 6\}}
\end{exbox}

\fname{ECPAnalyticRank}
\usageline{ECPAnalyticRank[model]\quad ECPAnalyticRank[model, prec]}
\args\ \texttt{model}: an integral long model $\{a_1,\dots,a_6\}$ (short pairs, equations and \texttt{Automatic} are also accepted). \texttt{prec} (default 40): working precision.
\returns\ An Association with keys \texttt{"Rank"} (order of vanishing of $L(E,s)$ at $s=1$), \texttt{"Parity"} (the root number $\pm1$) and \texttt{"Conductor"}.
\notes\ The rank is the true Mordell--Weil rank only assuming the Birch--Swinnerton-Dyer conjecture (proved for analytic rank $\le 1$). By the parity conjecture the root number alone decides whether the rank is even or odd. Precision below the internal minimum raises a message. The slowest command: its computation may take time (minutes for conductor $N\gtrsim10^5$), since it needs of order $\sqrt N$ coefficients $a_n$ and many high-precision incomplete gamma functions. A lower \texttt{prec} speeds it up.

\begin{exbox}
\textbf{Example.} analytic rank, root number and conductor of 37a1

\MatIn{6}{ECPAnalyticRank[\{0, 0, 1, -1, 0\}]}

\smallskip

\MatOut{6}{<|"Rank" -> 1, "Parity" -> -1, "Conductor" -> 37|>}
\end{exbox}

\fname{ECPAreProjectivelyEquivalent}
\usageline{ECPAreProjectivelyEquivalent[P, Q]}
\args\ \texttt{P}, \texttt{Q}: length-3 lists of numbers or symbols, the projective coordinates $[x:y:z]$.
\returns\ \texttt{True} if $Q=\lambda P$ for some $\lambda\neq 0$, \texttt{False} otherwise.
\notes\ The zero vector is equivalent to nothing, not even itself. The test does not check that the points lie on any curve.

\begin{exbox}
\textbf{Example.} proportional projective points

\MatIn{7}{ECPAreProjectivelyEquivalent[\{1, 2, 3\}, \{2, 4, 6\}]}

\smallskip

\MatOut{7}{True}
\end{exbox}

\fname{ECPBachetDuplication}
\usageline{ECPBachetDuplication[x, curve]\quad ECPBachetDuplication[\{x, y\}, curve]
ECPBachetDuplication[Infinity, curve]}
\args\ \texttt{x}: a number or a symbol. \texttt{\{x,y\}}: a point of the curve (checked). \texttt{curve}: any curve form.
\returns\ For a bare $x$: $x(2P)=\dfrac{x^4-b_4x^2-2b_6x-b_8}{4x^3+b_2x^2+2b_4x+b_6}$ (for $y^2=x^3+Ax+B$: $\dfrac{x^4-2Ax^2-8Bx+A^2}{4(x^3+Ax+B)}$; for $y^2=x^3+c$: $\dfrac{x^4-8cx}{4(x^3+c)}$). For a point: the whole point $2P$.
\notes\ The result is \texttt{Infinity} when $2P=\mathcal O$ (for a numeric $x$ at a 2-torsion point). The whole-point form agrees with \texttt{ECPMultiply[2, P, curve]}.

\begin{exbox}
\textbf{Example.} the whole point 2P

\MatIn{8}{ECPBachetDuplication[\{-3, 9\}, \{-36, 0\}]}

\smallskip

\MatOut{8}{\{25/4, -35/8\}}
\end{exbox}

\fname{ECPCanonicalHeight}
\usageline{ECPCanonicalHeight[P, curve]\quad ECPCanonicalHeight[P, curve, maxIter]}
\args\ \texttt{P}: a rational point on \texttt{curve}, or \texttt{Infinity}. \texttt{curve}: any curve form. \texttt{maxIter} (default 8): number of doublings in the limit $\hat h(P)=\tfrac12\lim h(2^nP)/4^n$, where $h$ is the naive height of \texttt{ECPNaiveHeight}.
\returns\ A machine-precision real number. Torsion points and \texttt{Infinity} give exactly $0$.
\notes\ The value is an approximation (only \texttt{maxIter} doublings), so different \texttt{maxIter} may disagree in the last digits; this is why \texttt{ECPIndependentQ} uses a loose default tolerance. \emph{Normalisation:} the package follows the convention of Silverman's \emph{Arithmetic of Elliptic Curves}, in which $\hat h$ is half of the value tabulated in Silverman--Tate, Cremona and the LMFDB. Its computation may take time for points of large height, since each doubling roughly quadruples the number of digits in the coordinates.

\begin{exbox}
\textbf{Example.} the generator (0,0) of 37a1 (half of Cremona's value $0.0511$)

\MatIn{9}{ECPCanonicalHeight[\{0, 0\}, \{0, 0, 1, -1, 0\}]}

\smallskip

\MatOut{9}{0.0255553}
\end{exbox}

\fname{ECPChordAndTangent}
\usageline{ECPChordAndTangent[P, curve]\quad ECPChordAndTangent[P, curve, vars]}
\args\ \texttt{P}: projective point $\{x,y,z\}$ on the cubic. \texttt{curve}: any curve form (homogenised internally). \texttt{vars}: the three variable names (default: inferred).
\returns\ A projective point: the third point where the tangent at $P$ meets the cubic, before any reflection. It is the chord-and-tangent doubling, not $2P$ in the group law.
\notes\ If $P$ is a flex the function returns $P$ itself; at a 2-torsion point the tangent is vertical and the result is the point at infinity. The computation is an exact projective parametrisation and does not depend on the chart.

\begin{exbox}
\textbf{Example.} third intersection of the tangent at (0:0:1)

\MatIn{10}{ECPChordAndTangent[\{0, 0, 1\}, y\textasciicircum{}2 z + y z\textasciicircum{}2 - x\textasciicircum{}3 + x z\textasciicircum{}2, \{x, y, z\}]}

\smallskip

\MatOut{10}{\{1, -1, 1\}}
\end{exbox}

\fname{ECPComposeTransformations}
\usageline{ECPComposeTransformations[\{u1,r1,s1,t1\}, \{u2,r2,s2,t2\}]}
\args\ Two isomorphisms, each a length-4 list with $u\neq 0$, for $x=u^2x'+r$, $y=u^3y'+u^2sx'+t$.
\returns\ One list $\{u,r,s,t\}$ equal to applying the first and then the second transformation.
\notes\ Composition is not commutative; the order of the arguments is the order of application. $u=0$ raises a message and returns \texttt{\$Failed}.

\begin{exbox}
\textbf{Example.} compose two isomorphisms

\MatIn{11}{ECPComposeTransformations[\{2, 1, 3, 4\}, \{5, 6, 7, 8\}]}

\smallskip

\MatOut{11}{\{10, 25, 17, 140\}}
\end{exbox}

\fname{ECPConductor}
\usageline{ECPConductor[model]}
\args\ \texttt{model}: an integral long model $\{a_1,\dots,a_6\}$.
\returns\ The integer $N=\prod_p p^{f_p}$ over the conductor exponents $f_p$ found by Tate's algorithm.
\notes\ Built on \texttt{ECPLocalData}. A non-integral model raises a message and returns \texttt{\$Failed}. Its computation may take time, since it factorises $\Delta$, which is slow when $|\Delta|$ has large prime factors.

\begin{exbox}
\textbf{Example.} conductor of 37a1

\MatIn{12}{ECPConductor[\{0, 0, 1, -1, 0\}]}

\smallskip

\MatOut{12}{37}
\end{exbox}

\fname{ECPCoulombDimension}
\usageline{ECPCoulombDimension[ordDelta]}
\args\ \texttt{ordDelta}: an integer with $1<\operatorname{ord}\Delta<12$.
\returns\ The scaling dimension $12/(12-\operatorname{ord}\Delta)$ of the Coulomb-branch parameter at a Kodaira singularity with that order of vanishing of $\Delta$ (Argyres--Douglas and Minahan--Nemeschansky: $II\to 6/5$, $III\to 4/3$, $IV\to 3/2$, $I_0^*\to 2$, $IV^*\to 3$, $III^*\to 4$, $II^*\to 6$).
\notes\ Meaningful only for the types that are not $I_n$ or $I_n^*$, that is, for $\operatorname{ord}\Delta\in\{2,3,4,6,8,9,10\}$.

\begin{exbox}
\textbf{Example.} the seven values

\MatIn{93}{ECPCoulombDimension[\#] \& /@ \{2, 3, 4, 6, 8, 9, 10\}}

\smallskip

\MatOut{93}{\{6/5, 4/3, 3/2, 2, 3, 4, 6\}}
\end{exbox}

\fname{ECPCubicToWeierstrass}
\usageline{ECPCubicToWeierstrass[P, curve]\quad ECPCubicToWeierstrass[P, curve, opts]}
\args\ \texttt{P}: a rational point $\{x_0,y_0,z_0\}$ of the cubic. \texttt{curve}: a homogeneous cubic in three variables, or an affine cubic (homogenised with respect to $z$). \texttt{opts}: \texttt{"Form"} (\texttt{"Long"}, default, or \texttt{"Short"}) and the algorithm/equivalence options.
\returns\ The long Weierstrass polynomial $y^2+a_1xy+a_3y-(x^3+a_2x^2+a_4x+a_6)$; with \texttt{"Form" -> "Short"} the equation $Y^2=X^3+AX+B$. The option \texttt{"Equivalence" -> True} also prints the transformation used.
\notes\ $P$ need not be a flex. For characteristic 2 or 3 use \texttt{ECPCubicToWeierstrassCharTwoThree}. A point that is not on the cubic, or an invalid option, raises a message and returns \texttt{\$Failed}.

\begin{exbox}
\textbf{Example.} cubic with the rational point (0:0:1)

\MatIn{13}{ECPCubicToWeierstrass[\{0, 0, 1\}, y\textasciicircum{}2 z + y z\textasciicircum{}2 - x\textasciicircum{}3 + x z\textasciicircum{}2]}

\smallskip

\MatOut{13}{7 x\textasciicircum{}2 - x\textasciicircum{}3 + y + 4 x y + y\textasciicircum{}2}
\end{exbox}

\fname{ECPCubicToWeierstrassCharTwoThree}
\usageline{ECPCubicToWeierstrassCharTwoThree[F, P, p]\quad ECPCubicToWeierstrassCharTwoThree[F, P, \{p, r\}]}
\args\ \texttt{F}: a homogeneous cubic in $x,y,z$. \texttt{P}: a flex of \texttt{F} (the point $[0:1:0]$ of a Weierstrass cubic always works). \texttt{p}: 2 or 3; \texttt{\{p,r\}}: work over $\mathbb F_{p^r}$.
\returns\ An Association with \texttt{"WeierstrassCoefficients"}, \texttt{"Transformation"} and auxiliary data (over $\mathbb F_{p^r}$ also \texttt{"Scaling"}).
\notes\ A purely linear change of coordinates anchored at the supplied flex; the usual completion of the square or cube would divide by 2 or 3.

\begin{exbox}
\textbf{Example.} cubic over F\_3 with the flex (0:1:0)

\MatIn{14}{ECPCubicToWeierstrassCharTwoThree[y\textasciicircum{}2 z - x\textasciicircum{}3 - x z\textasciicircum{}2 - z\textasciicircum{}3, \{0, 1, 0\}, 3]}

\smallskip

\MatOut{14}{<|"WeierstrassCoefficients" -> \{0, 0, 0, 1, 1\}, "Transformation" -> \{\{1, 0, 0\}, \{0, 1, 0\}, \{0, 0, 1\}\}, "AuxiliaryPoints" -> <|"Q" -> \{1, 0, 0\}, "P" -> \{0, 1, 0\}, "R" -> \{0, 0, 1\}|>, "Characteristic" -> 3|>}
\end{exbox}

\fname{ECPElGamalDecodeFp}
\usageline{ECPElGamalDecodeFp[s, \{C1,C2\}, a, p]\quad ECPElGamalDecodeFp[s, \{C1,C2\}, curve, p, kappa]}
\args\ \texttt{s}: the private key. \texttt{\{C1,C2\}}: the ciphertext from \texttt{ECPElGamalEncodeFp}. \texttt{a} or \texttt{curve}: the curve. \texttt{p}: odd prime. \texttt{kappa} (default 100): the Koblitz parameter used at encoding time.
\returns\ The plaintext integer $m$.
\notes\ Computes $M=C_2-sC_1$ and inverts the embedding $x=m\kappa+j$. A wrong \texttt{kappa} silently recovers a wrong integer. A field $\{p,r\}$ with $r\ge2$ raises a message.

\begin{exbox}
\textbf{Example.} decrypt with the secret key s = 57 (round trip returns 42)

\MatIn{15}{ECPElGamalDecodeFp[57, ECPElGamalEncodeFp[42, \{\{0, 8530\}, \{1753, 557\}\}, 2, 3, 10007, 4923], 2, 10007]}

\smallskip

\MatOut{15}{42}
\end{exbox}

\fname{ECPElGamalEncodeFp}
\usageline{ECPElGamalEncodeFp[m, \{P,B\}, a, b, p, n]\quad ECPElGamalEncodeFp[m, \{P,B\}, curve, p, n, kappa]}
\args\ \texttt{m}: non-negative integer. \texttt{\{P,B\}}: base point and public key $B=sP$. \texttt{a,b} or \texttt{curve}: the curve. \texttt{p}: odd prime. \texttt{n}: bound for the ephemeral $k$ (normally the order of $P$). \texttt{kappa} (default 100): Koblitz parameter.
\returns\ The ciphertext $\{C_1,C_2\}$ with $C_1=kP$, $C_2=M+kB$, where $M$ is the Koblitz point embedding $m$.
\notes\ Requires $\kappa(m+1)<p$. The ephemeral key $k$ is drawn at random from $[\lfloor n/4\rfloor,n-1]$, so the ciphertext changes from one call to the next. Toy primes only in this manual.

\begin{exbox}
\textbf{Example.} encrypt $m=42$ on $y^2=x^3+2x+3$ over $\mathbb F_{10007}$ with $P=(0,8530)$ of order $4923$ and $B=57P=(1753,557)$; the ciphertext depends on the random $k$, so we check it by decrypting

\MatIn{16}{c = ECPElGamalEncodeFp[42, \{\{0, 8530\}, \{1753, 557\}\}, 2, 3, 10007, 4923]; ECPElGamalDecodeFp[57, c, 2, 10007]}

\smallskip

\MatOut{16}{42}
\end{exbox}

\fname{ECPExtractCoefficients}
\usageline{ECPExtractCoefficients[curve]}
\args\ \texttt{curve}: a polynomial, an equation, a coefficient list or an Association; variables $x,y$ or $X,Y$.
\returns\ The five long coefficients.
\notes\ A singular curve raises a message and returns \texttt{\$Failed}; this is the check every other command relies on.

\begin{exbox}
\textbf{Example.} read the coefficients off an equation

\MatIn{17}{ECPExtractCoefficients[y\textasciicircum{}2 + y == x\textasciicircum{}3 - x]}

\smallskip

\MatOut{17}{\{0, 0, 1, -1, 0\}}
\end{exbox}

\fname{ECPFibreAt}
\usageline{ECPFibreAt[curve, t, t0]}
\args\ \texttt{curve}: a family over $\mathbb Q(t)$ in any accepted form (long $\{a_1,\dots,a_6\}$ with coefficients that are polynomials in \texttt{t}, short $\{A,B\}$, an equation in $x,y$, an Association). \texttt{t}: the parameter (a symbol). \texttt{t0}: a number, or \texttt{Infinity} (then \texttt{ECPFibreAtInfinity} is used).
\returns\ An Association with the keys \texttt{"Location"}, \texttt{"Degree"}, \texttt{"Type"} (Kodaira symbol), \texttt{"Euler"}, \texttt{"Orders"} ($\{\operatorname{ord}c_4,\operatorname{ord}c_6,\operatorname{ord}\Delta\}$), \texttt{"Split"}, \texttt{"SplitCriterion"}, \texttt{"GaugeAlgebra"} and \texttt{"ScalingDimension"}.
\notes\ The model must be minimal at \texttt{t0}; otherwise the type is reported as \texttt{"non-minimal"}. \texttt{"Split"} is decided over $\mathbb Q$ for the types $I_n$ ($n\ge1$) and $IV$ by testing whether \texttt{"SplitCriterion"} ($-c_6$ at the point for $I_n$, the leading coefficient of $-c_6/864$ for $IV$) is a rational square; for symbolic models it is \texttt{Missing["Symbolic"]} and \texttt{"SplitCriterion"} is the expression that must be a square. For other types it is \texttt{Missing["NotApplicable"]}. The invariants $c_4,c_6,\Delta$ must be polynomials in \texttt{t}; a singular family raises a message.

\begin{exbox}
\textbf{Example.} the $E_8$ Tate form has a fibre of type $II^*$ at $t=0$

\MatIn{89}{Lookup[ECPFibreAt[ECPTateForm[\{1, 2, 3, 4, 5\}, t], t, 0], \{"Type", "Euler", "GaugeAlgebra", "ScalingDimension"\}]}

\smallskip

\MatOut{89}{\{"II*", 10, "e8", 6\}}
\end{exbox}

\fname{ECPFibreAtInfinity}
\usageline{ECPFibreAtInfinity[curve, t]\quad ECPFibreAtInfinity[curve, t, n]}
\args\ \texttt{curve}, \texttt{t} as for \texttt{ECPFibreAt}; \texttt{n}: a positive integer or \texttt{Automatic} (default).
\returns\ The same Association as \texttt{ECPFibreAt}, for the fibre at $t=\infty$ of a family over $\mathbb P^1$ whose coefficients satisfy $\deg a_i\le in$ ($n=1$: rational elliptic surface, $n=2$: K3).
\notes\ The fibre is classified in $s=1/t$ with the weights $s^{4n}c_4$, $s^{6n}c_6$, $s^{12n}\Delta$. With \texttt{Automatic}, $n$ is the smallest value that makes the model regular at infinity; give $n$ explicitly to use a larger one. A value of $n$ that is too small raises a message.

\begin{exbox}
\textbf{Example.} the universal curve with a point of order 6 has an $I_2$ fibre at infinity

\MatIn{90}{Lookup[ECPFibreAtInfinity[ECPTorsionFamily[6, d]["Coefficients"], d], \{"Type", "Euler"\}]}

\smallskip

\MatOut{90}{\{"I2", 2\}}
\end{exbox}

\fname{ECPFibreList}
\usageline{ECPFibreList[curve, t]\quad ECPFibreList[curve, t, n]}
\args\ \texttt{curve}, \texttt{t}, \texttt{n} as for \texttt{ECPFibreAtInfinity}.
\returns\ An Association with the keys \texttt{"n"}, \texttt{"Fibres"} (one record as in \texttt{ECPFibreAt} for every irreducible factor of $\Delta$ over $\mathbb Q$, with \texttt{"Degree"} its number of points, followed by the fibre at \texttt{Infinity}), \texttt{"TotalEuler"}, \texttt{"ExpectedEuler"} ($=12n$) and \texttt{"ConsistentQ"}.
\notes\ \texttt{"ConsistentQ"} is \texttt{True} when the Euler numbers add up to $12n$, a strong check that the model is minimal and the classification is right. For an irreducible factor of degree larger than one the \texttt{"Location"} is that polynomial (its roots all carry the same fibre type) and \texttt{"Split"} is \texttt{Missing["NeedsExtension"]}.

\begin{exbox}
\textbf{Example.} the family with a point of order 6 has fibres $I_6,I_3,I_2,I_1$

\MatIn{91}{Sort[Lookup[\#, \{"Type", "Location"\}] \& /@ ECPFibreList[ECPTorsionFamily[6, d]["Coefficients"], d]["Fibres"]]}

\smallskip

\MatOut{91}{\{\{"I1", -1/9\}, \{"I2", Infinity\}, \{"I3", -1\}, \{"I6", 0\}\}}

\medskip

\MatIn{92}{Lookup[ECPFibreList[ECPTorsionFamily[6, d]["Coefficients"], d], \{"n", "TotalEuler", "ExpectedEuler", "ConsistentQ"\}]}

\smallskip

\MatOut{92}{\{1, 12, 12, True\}}
\end{exbox}

\fname{ECPFindGenerators}
\usageline{ECPFindGenerators[\{A,B\}]\quad ECPFindGenerators[curve, opts]}
\args\ \texttt{curve}: any curve form; the search runs on the integral short model and results are translated back. \texttt{opts}: \texttt{"SearchBound"} (default 50), \texttt{"Verbose"} (default \texttt{False}).
\returns\ An Association with \texttt{"Rank"} (the number of independent points found, a lower bound for the rank), \texttt{"Generators"}, \texttt{"Torsion"}, \texttt{"RankBound"} and diagnostic data. The key \texttt{"RankBound"} is the heuristic estimate of \texttt{ECPRankUpperBound}, not a proven bound; for a proven bound use \texttt{ECPRankBound} (which gives $3$ for the example below). 
\notes\ A bounded heuristic, not a 2-descent: generators outside the search range are missed and the rank is not proven. The generators returned generate the free part of the Mordell--Weil group (modulo torsion), not the torsion subgroup. Its computation may take time, since a canonical height is computed for every point found; the cost grows with \texttt{"SearchBound"}. In the example below, the generator $(-2,-8)=(-3,9)+(6,0)$ differs from the point $(-3,9)$ used elsewhere in this manual by a torsion point, so both have the same canonical height.

\begin{exbox}
\textbf{Example.} independent generators, short form

\MatIn{18}{ECPFindGenerators[\{-36, 0\}, "SearchBound" -> 30]}

\smallskip

\MatOut{18}{<|"Rank" -> 1, "Generators" -> \{\{-2, -8\}\}, "Torsion" -> \{Infinity, \{-6, 0\}, \{0, 0\}, \{6, 0\}\}, "RankBound" -> 2, "Curve" -> \{-36, 0\}, "HeightMatrix" -> \{\{0.444313\}\}, "Determinant" -> 0.444313, "Time" -> 62.1973|>}
\end{exbox}

\fname{ECPFindLocalReduction}
\usageline{ECPFindLocalReduction[\{a1,...,a6\}, p]}
\args\ A long model with integer coefficients and a prime $p$.
\returns\ A list $\{u,r,s,t\}$ that reduces the model at $p$ by one step; the identity $\{1,0,0,0\}$ if the model is already minimal at $p$.
\notes\ One step at one prime; \texttt{ECPMinimalModel} iterates it over all primes.

\begin{exbox}
\textbf{Example.} Laska-Kraus step at the prime 2

\MatIn{19}{ECPFindLocalReduction[\{0, 0, 0, -576, 0\}, 2]}

\smallskip

\MatOut{19}{\{2, 0, 0, 0\}}
\end{exbox}

\fname{ECPFindPointFp}
\usageline{ECPFindPointFp[m, a, b, p]\quad ECPFindPointFp[m, curve, p, kappa]}
\args\ \texttt{m}: non-negative integer message. \texttt{a,b} or \texttt{curve}: the curve over an odd prime field. \texttt{p}: odd prime. \texttt{kappa} (default 100).
\returns\ A point $\{x,y\}$ with $x=m\kappa+j$ for the smallest successful $j$.
\notes\ Requires $\kappa(m+1)<p$. $\{p,r\}$ with $r\ge2$ raises a message.

\begin{exbox}
\textbf{Example.} embed the message 3 as a point (Koblitz)

\MatIn{20}{ECPFindPointFp[3, 2, 3, 10007]}

\smallskip

\MatOut{20}{\{300, 349\}}
\end{exbox}

\fname{ECPFindPointOrder}
\usageline{ECPFindPointOrder[point, a, b]\quad ECPFindPointOrder[point, curve, maxOrder]}
\args\ \texttt{point}: $\{x,y\}$ or \texttt{Infinity}. \texttt{a,b} or \texttt{curve}. \texttt{maxOrder} (default 12).
\returns\ The smallest $n\le$\texttt{maxOrder} with $nP=\mathcal O$, otherwise \texttt{Infinity}.
\notes\ Brute force; compare \texttt{ECPPointOrderFp}, which is exact over $\mathbb F_p$.

\begin{exbox}
\textbf{Example.} order of the 2-torsion point (0,0), short coefficients a, b

\MatIn{21}{ECPFindPointOrder[\{0, 0\}, -36, 0]}

\smallskip

\MatOut{21}{2}
\end{exbox}

\fname{ECPFindPointOrderFp}
\usageline{ECPFindPointOrderFp[point, curve, p]\quad ECPFindPointOrderFp[point, curve, \{p, r\}]}
\args\ A point or \texttt{Infinity} on the curve mod $p$; any curve form (a short pair needs $p$ odd); a prime $p$ or $\{p,r\}$.
\returns\ The exact order of the point: the smallest divisor of $N=\#E(\mathbb F_q)$ that kills it.
\notes\ Same algorithm as \texttt{ECPPointOrderFp} (the group order comes from \texttt{ECPSEAPointCount}, so about $p\le10^6$, or $q\le$ \texttt{\$ECPFqTableLimit}). Returns \texttt{\$Failed} with a message for a non-prime, a singular curve or a point off the curve.

\begin{exbox}
\textbf{Example.} order of (5,1) on $y^2=x^3+2x+2$ over $\mathbb F_{17}$

\MatIn{98}{ECPFindPointOrderFp[\{5, 1\}, \{2, 2\}, 17]}

\smallskip

\MatOut{98}{19}
\end{exbox}

\fname{ECPFindTorsionPoints}
\usageline{ECPFindTorsionPoints[curve]}
\args\ \texttt{curve}: any curve form.
\returns\ A list of all rational torsion points, \texttt{Infinity} first, in the coordinates of \texttt{curve}.
\notes\ Nagell--Lutz: a nonzero torsion point has integral coordinates with $y=0$ or $y^2\mid\Delta$. For non-integral coefficients or $|\Delta|>$\texttt{\$ECPNagellLutzDiscLimit} the function falls back to a bounded search and says so; the result is then not guaranteed complete. Its computation may take time when $|\Delta|$ is large, since $\Delta$ is factorised and all its divisors are listed.

\begin{exbox}
\textbf{Example.} Nagell-Lutz on y\textasciicircum{}2 = x\textasciicircum{}3 - 36x

\MatIn{22}{ECPFindTorsionPoints[\{-36, 0\}]}

\smallskip

\MatOut{22}{\{Infinity, \{-6, 0\}, \{0, 0\}, \{6, 0\}\}}
\end{exbox}

\fname{ECPFiniteFieldInfo}
\usageline{ECPFiniteFieldInfo[\{p, r\}]}
\args\ \texttt{\{p,r\}}: a prime and $r\ge2$.
\returns\ The field polynomial $f(x)$ over $\mathbb F_p$ (so $\mathbb F_{p^r}\cong\mathbb F_p[x]/(f)$) and the meaning of the generator $g$ (the class of $x$, a primitive element).
\notes\ Elements are written $0$ or $g^k$ ($g^0$ prints as $1$). On input any Laurent polynomial in $g$ with integer coefficients is accepted (\texttt{1+g}, \texttt{2 g\^{}3+1}, \texttt{g\^{}-1}), an integer $n$ means $n$ mod $p$, and any single symbol may replace $g$.

\begin{exbox}
\textbf{Example.} how F\_8 is represented

\MatIn{23}{ECPFiniteFieldInfo[\{2, 3\}]}

\smallskip

\MatOut{23}{<|"p" -> 2, "r" -> 3, "q" -> 8, "FieldPolynomial" -> 1 + x + x\textasciicircum{}3, "Generator" -> "g, the class of x modulo FieldPolynomial; g is a primitive element (it generates F\_q\textasciicircum{}*)", "Elements" -> "0 and g\textasciicircum{}k for k = 0..q-2 (g\textasciicircum{}0 prints as 1)"|>}
\end{exbox}

\fname{ECPFormattedEquation}
\usageline{ECPFormattedEquation[curve]}
\args\ \texttt{curve}: long $\{a_1,\dots,a_6\}$, short $\{A,B\}$, Association, equation, polynomial or \texttt{Automatic}.
\returns\ The equation $y^2+a_1xy+a_3y=x^3+a_2x^2+a_4x+a_6$ as an \texttt{Equal} expression.
\notes\ Display only; use \texttt{ECPExtractCoefficients} to go the other way.

\begin{exbox}
\textbf{Example.} display the long Weierstrass equation

\MatIn{24}{ECPFormattedEquation[\{0, 0, 1, -1, 0\}]}

\smallskip

\MatOut{24}{y + y\textasciicircum{}2 == -x + x\textasciicircum{}3}
\end{exbox}

\fname{ECPFqAdd}
\usageline{ECPFqAdd[a, b, \{p, r\}]}
\args\ Field elements in the $0/g^k$ notation (or polynomials in $g$); $r\ge2$.
\returns\ The sum, as $0$ or $g^k$.
\notes\ Uses the discrete-log tables, so $q=p^r\le$ \texttt{\$ECPFqTableLimit}.

\begin{exbox}
\textbf{Example.} g + g\textasciicircum{}2 in F\_8

\MatIn{25}{ECPFqAdd[g, g\textasciicircum{}2, \{2, 3\}]}

\smallskip

\MatOut{25}{g\textasciicircum{}4}
\end{exbox}

\fname{ECPFqElement}
\usageline{ECPFqElement[\{c0, ..., c\_(r-1)\}, \{p, r\}]}
\args\ The coefficients $c_i\in\mathbb F_p$ and the field $\{p,r\}$.
\returns\ The element $c_0+c_1g+\dots+c_{r-1}g^{r-1}$ in $0/g^k$ notation.
\notes\ Inverse of \texttt{ECPFqVector}.

\begin{exbox}
\textbf{Example.} the element with coordinates {1,1,0}

\MatIn{26}{ECPFqElement[\{1, 1, 0\}, \{2, 3\}]}

\smallskip

\MatOut{26}{g\textasciicircum{}3}
\end{exbox}

\fname{ECPFqInverse}
\usageline{ECPFqInverse[a, \{p, r\}]}
\args\ A nonzero element $a$ and the field $\{p,r\}$.
\returns\ $a^{-1}$ as $g^k$.
\notes\ Zero has no inverse and raises a message.

\begin{exbox}
\textbf{Example.} inverse of g in F\_8

\MatIn{27}{ECPFqInverse[g, \{2, 3\}]}

\smallskip

\MatOut{27}{g\textasciicircum{}6}
\end{exbox}

\fname{ECPFqMultiply}
\usageline{ECPFqMultiply[a, b, \{p, r\}]}
\args\ Field elements and the field $\{p,r\}$.
\returns\ The product, as $0$ or $g^k$.
\notes\ Multiplication is addition of exponents modulo $q-1$.

\begin{exbox}
\textbf{Example.} g\textasciicircum{}3 g\textasciicircum{}5 in F\_8

\MatIn{28}{ECPFqMultiply[g\textasciicircum{}3, g\textasciicircum{}5, \{2, 3\}]}

\smallskip

\MatOut{28}{g}
\end{exbox}

\fname{ECPFqVector}
\usageline{ECPFqVector[a, \{p, r\}]}
\args\ A field element and the field $\{p,r\}$.
\returns\ The list $\{c_0,\dots,c_{r-1}\}$ over $\mathbb F_p$ with $a=\sum c_ig^i$.
\notes\ Inverse of \texttt{ECPFqElement}.

\begin{exbox}
\textbf{Example.} coordinates of g\textasciicircum{}3 over F\_2

\MatIn{29}{ECPFqVector[g\textasciicircum{}3, \{2, 3\}]}

\smallskip

\MatOut{29}{\{1, 1, 0\}}
\end{exbox}

\fname{ECPFromMinimalPoint}
\usageline{ECPFromMinimalPoint[\{X,Y\}, \{u,r,s,t\}]}
\args\ A point on the transformed model and the isomorphism that produced it (e.g.\ the \texttt{"Transformation"} of \texttt{ECPMinimalModel}).
\returns\ The point $\{x,y\}$ on the original model: $x=u^2X+r$, $y=u^3Y+u^2sX+t$.
\notes\ Exact inverse of \texttt{ECPToMinimalPoint}; $u=0$ raises a message.

\begin{exbox}
\textbf{Example.} and back again

\MatIn{30}{ECPFromMinimalPoint[ECPToMinimalPoint[\{-12, 72\}, ECPMinimalModel[\{0, 0, 0, -576, 0\}]["Transformation"]], ECPMinimalModel[\{0, 0, 0, -576, 0\}]["Transformation"]]}

\smallskip

\MatOut{30}{\{-12, 72\}}
\end{exbox}

\fname{ECPFromShortCoordinates}
\usageline{ECPFromShortCoordinates[\{X,Y\}, curve]\quad ECPFromShortCoordinates[Infinity, curve]}
\args\ A point on the short model $Y^2=X^3+AX+B$ of \texttt{ToShortWeierstrass[curve]}; \texttt{curve} is the original long curve.
\returns\ The point on \texttt{curve}: the inverse of $X=x+b_2/12$, $Y=y+(a_1x+a_3)/2$.
\notes\ Exact inverse of \texttt{ECPToShortCoordinates}. A point not on the short model raises a message.

\begin{exbox}
\textbf{Example.} the inverse map

\MatIn{31}{ECPFromShortCoordinates[\{0, 1/2\}, \{0, 0, 1, -1, 0\}]}

\smallskip

\MatOut{31}{\{0, 0\}}
\end{exbox}

\fname{ECPFromShortCoordinatesFp}
\usageline{ECPFromShortCoordinatesFp[\{X,Y\}, curve, p]}
\args\ A point on the model returned by \texttt{ECPToShortWeierstrassFp[curve, p]}; the original long curve; a prime $p$.
\returns\ The corresponding point of the original curve mod $p$.
\notes\ Works in every characteristic, unlike \texttt{ECPFromShortCoordinates}.

\begin{exbox}
\textbf{Example.} and back again

\MatIn{32}{ECPFromShortCoordinatesFp[ECPToShortCoordinatesFp[\{0, 0\}, \{0, 0, 1, -1, 0\}, 7], \{0, 0, 1, -1, 0\}, 7]}

\smallskip

\MatOut{32}{\{0, 0\}}
\end{exbox}

\fname{ECPGaugeAlgebra}
\usageline{ECPGaugeAlgebra[type]\quad ECPGaugeAlgebra[type, split]}
\args\ \texttt{type}: a Kodaira symbol as a string (\texttt{"I3"}, \texttt{"IV*"}, \dots). \texttt{split}: \texttt{True} or \texttt{False} where it matters (default \texttt{Automatic}).
\returns\ The F-theory gauge algebra as a string: $I_n$: su($n$) / sp($\lfloor n/2\rfloor$); $IV$: su(3) / su(2); $I_n^*$: so($2n+8$) / so($2n+7$); $IV^*$: e6 / f4; $III$ and $I_2$: su(2); $III^*$: e7; $II^*$: e8.
\notes\ $I_0^*$ depends on the monodromy cover (so(8), so(7) or $\mathfrak g_2$), and $II$, $I_0$, $I_1$ carry no gauge algebra. When \texttt{split} is omitted, an unresolved case is returned as an ``either/or'' string.

\begin{exbox}
\textbf{Example.} split and non-split $IV^*$, and an unresolved $I_5$

\MatIn{94}{\{ECPGaugeAlgebra["IV*", True], ECPGaugeAlgebra["IV*", False], ECPGaugeAlgebra["I5"]\}}

\smallskip

\MatOut{94}{\{"e6", "f4", "su(5) or sp(2)"\}}
\end{exbox}

\fname{ECPGroupFp}
\usageline{ECPGroupFp[curve, p]\quad ECPGroupFp[curve, \{p, r\}]\quad ECPGroupFp[a, b, p]}
\args\ \texttt{curve}: short, long, Association, equation or \texttt{Automatic}. \texttt{p} or \texttt{\{p,r\}}: the field, $q=p^r\le$ \texttt{\$ECPFqTableLimit}. \texttt{a,b}: the short curve $y^2=x^3+ax+b$.
\returns\ An Association with \texttt{"Field"}, \texttt{"Curve"}, \texttt{"Order"}, \texttt{"Invariants"} $\{n_1,n_2\}$ with $n_1\mid n_2$ and $E(\mathbb F_q)\cong\mathbb Z/n_1\times\mathbb Z/n_2$, \texttt{"Structure"} (a string), \texttt{"Generators"} (one point if cyclic, two otherwise) and \texttt{"Points"} (all points, \texttt{Infinity} first).
\notes\ Works in every characteristic. Over $\mathbb F_p$ points are integer pairs; over $\mathbb F_{p^r}$ coordinates are $0$ or $g^k$. Its computation may take time for large $q$, since every point of $E(\mathbb F_q)$ is listed.

\begin{exbox}
\textbf{Example.} a non-cyclic group over F\_7

\MatIn{33}{ECPGroupFp[\{-1, 0\}, 7]["Structure"]}

\smallskip

\MatOut{33}{"Z/2 x Z/4"}
\end{exbox}

\fname{ECPHeightMatrix}
\usageline{ECPHeightMatrix[points, curve]}
\args\ A list of points on \texttt{curve}; \texttt{curve}: any curve form.
\returns\ The symmetric matrix $M_{ij}=\langle P_i,P_j\rangle$.
\notes\ Torsion points give zero rows; dependent points give a determinant near $0$. Its computation may take time, since $k$ points need about $3k^2$ canonical heights.

\begin{exbox}
\textbf{Example.} height matrix; its determinant is the regulator

\MatIn{34}{ECPHeightMatrix[\{\{2, 9\}, \{3, 10\}\}, \{0, 73\}]}

\smallskip

\MatOut{34}{\{\{0.46197, -0.488523\}, \{-0.488523, 0.996347\}\}}
\end{exbox}

\fname{ECPHeightPairing}
\usageline{ECPHeightPairing[P1, P2, curve]}
\args\ Two points on \texttt{curve}; \texttt{curve}: any curve form.
\returns\ The real number $\langle P_1,P_2\rangle=\tfrac12\bigl(\hat h(P_1+P_2)-\hat h(P_1)-\hat h(P_2)\bigr)$, so that $\langle P,P\rangle=\hat h(P)$ exactly as returned by \texttt{ECPCanonicalHeight}.
\notes\ Inherits the approximation and the normalisation of \texttt{ECPCanonicalHeight}. Its computation may take time, since it needs three canonical heights.

\begin{exbox}
\textbf{Example.} pairing of (2,9) and (3,10) on y\textasciicircum{}2 = x\textasciicircum{}3 + 73

\MatIn{35}{ECPHeightPairing[\{2, 9\}, \{3, 10\}, \{0, 73\}]}

\smallskip

\MatOut{35}{-0.488523}
\end{exbox}

\fname{ECPIndependentQ}
\usageline{ECPIndependentQ[points, curve]\quad ECPIndependentQ[points, curve, tol]}
\args\ A list of points; \texttt{curve}; \texttt{tol} (default 0.01): determinant threshold.
\returns\ \texttt{True} if $\det$\,\texttt{ECPHeightMatrix} $>$ \texttt{tol}, else \texttt{False}.
\notes\ Any torsion point, or two multiples of one point, makes the set dependent. A negative \texttt{tol} raises a message. Its computation may take time, since it builds the full height matrix.

\begin{exbox}
\textbf{Example.} (2,9), (3,10) are independent

\MatIn{36}{ECPIndependentQ[\{\{2, 9\}, \{3, 10\}\}, \{0, 73\}]}

\smallskip

\MatOut{36}{True}
\end{exbox}

\fname{ECPIsMinimal}
\usageline{ECPIsMinimal[\{a1,...,a6\}]\quad ECPIsMinimal[curve]}
\args\ Any curve form (a short pair $\{A,B\}$ is read as $\{0,0,0,A,B\}$).
\returns\ \texttt{True} iff the long model is a global minimal (Laska--Kraus) model over $\mathbb Z$.
\notes\ The right test for a long model; \texttt{ECPMinimalQ} on a short pair can be fooled at $p=2,3$.

\begin{exbox}
\textbf{Example.} 37a1 has squarefree discriminant, so it is minimal

\MatIn{37}{ECPIsMinimal[\{0, 0, 1, -1, 0\}]}

\smallskip

\MatOut{37}{True}
\end{exbox}

\fname{ECPIsMinimalShort}
\usageline{ECPIsMinimalShort[\{A,B\}]\quad ECPIsMinimalShort[curve]}
\args\ A short pair, or any curve form.
\returns\ For a short pair: minimality among short models (via \texttt{ECPMinimalQ}); for long input: the same as \texttt{ECPIsMinimal}.
\notes\ A short pair can pass this test and still be non-minimal as a Weierstrass model at $p=2,3$.

\begin{exbox}
\textbf{Example.} y\textasciicircum{}2 = x\textasciicircum{}3 - 576x is not minimal (u = 2 removes 2\textasciicircum{}4, 2\textasciicircum{}6)

\MatIn{38}{ECPIsMinimalShort[\{-576, 0\}]}

\smallskip

\MatOut{38}{False}
\end{exbox}

\fname{ECPJFunction}
\usageline{ECPJFunction[curve]}
\args\ \texttt{curve}: a family in any accepted form, with coefficients that may contain symbols.
\returns\ The $j$-invariant $c_4^3/\Delta$ as a simplified rational function of the parameters.
\notes\ \texttt{ECPJInvariant} is meant for numeric curves; its singularity test cannot decide symbolic discriminants, which is why this command exists. A family with $\Delta=0$ raises a message.

\begin{exbox}
\textbf{Example.} the $j$-function of the Legendre family

\MatIn{95}{Simplify[ECPJFunction[\{0, -(1 + t), 0, t, 0\}] == 256 (t\textasciicircum{}2 - t + 1)\textasciicircum{}3/(t\textasciicircum{}2 (t - 1)\textasciicircum{}2)]}

\smallskip

\MatOut{95}{True}
\end{exbox}

\fname{ECPJInvariant}
\usageline{ECPJInvariant[curve]\quad ECPJInvariant[curve, p]\quad ECPJInvariant[curve, \{p, r\}]}
\args\ \texttt{curve}: any curve form (also bare coefficients). \texttt{p}: optional prime modulus; \texttt{\{p,r\}}: $\mathbb F_{p^r}$.
\returns\ $j=c_4^3/\Delta$, reduced mod $p$ or computed in $\mathbb F_{p^r}$ if a modulus is given.
\notes\ The short-form overload modulo $p$ excludes $p=2$.

\begin{exbox}
\textbf{Example.} curve = Automatic uses the session curve

\MatIn{39}{ECPJInvariant[Automatic]}

\smallskip

\MatOut{39}{110592/37}
\end{exbox}

\fname{ECPKodairaType}
\usageline{ECPKodairaType[\{ordC4, ordC6, ordDelta\}]}
\args\ The orders of vanishing of $c_4$, $c_6$ and $\Delta$ at a point: non-negative integers (or \texttt{Infinity}).
\returns\ \texttt{\{type, EulerNumber\}} from Kodaira's table (Table~\ref{tab:kodaira}); the Euler number equals $\operatorname{ord}\Delta$.
\notes\ The table is the characteristic-0 one, so it is valid over $\mathbb Q(t)$ and $\mathbb C(t)$. The result is \texttt{\{"non-minimal", Indeterminate\}} when $\operatorname{ord}c_4\ge4$ and $\operatorname{ord}c_6\ge6$, and \texttt{\{"unclassified", Indeterminate\}} for inconsistent orders.

\begin{exbox}
\textbf{Example.}

\MatIn{88}{\{ECPKodairaType[\{2, 3, 8\}], ECPKodairaType[\{4, 6, 12\}]\}}

\smallskip

\MatOut{88}{\{\{"I2*", 8\}, \{"non-minimal", Indeterminate\}\}}
\end{exbox}

\fname{ECPLegendreToWeierstrass}
\usageline{ECPLegendreToWeierstrass[lambda]}
\args\ \texttt{lambda}: a number or symbol, $\lambda\neq0,1$.
\returns\ A pair $\{A,B\}$ such that $y^2=x^3+Ax+B$ is isomorphic to $y^2=x(x-1)(x-\lambda)$ (scaling $A=-27c_4$, $B=-54c_6$).
\notes\ Requires characteristic 0 or $>3$.

\begin{exbox}
\textbf{Example.} Legendre form y\textasciicircum{}2 = x(x-1)(x-3)

\MatIn{40}{ECPLegendreToWeierstrass[3]}

\smallskip

\MatOut{40}{\{-3024, -34560\}}
\end{exbox}

\fname{ECPLocalData}
\usageline{ECPLocalData[model]}
\args\ An integral long model.
\returns\ An Association mapping each prime dividing $\Delta$ to the Association that \texttt{ECPTate[model, p]} returns.
\notes\ Equivalent to calling \texttt{ECPTate} once per prime factor of $\Delta$. Its computation may take time, since it factorises $\Delta$, which is slow when $|\Delta|$ has large prime factors.

\begin{exbox}
\textbf{Example.} reduction data at every bad prime

\MatIn{41}{ECPLocalData[\{0, 0, 0, -36, 0\}]}

\smallskip

\MatOut{41}{<|2 -> <|"f" -> 6, "type" -> "I2*", "model" -> \{0, 6, 0, -24, -64\}, "split" -> None|>, 3 -> <|"f" -> 2, "type" -> "I0*", "model" -> \{0, 0, 0, -36, 0\}, "split" -> None|>|>}
\end{exbox}

\fname{ECPMillerFp}
\usageline{ECPMillerFp[n, P, Q, curve, field]}
\args\ \texttt{n}: a positive integer. \texttt{P}, \texttt{Q}: points on \texttt{curve} over the field, with $Q\neq\mathcal O$. \texttt{curve}: short $\{A,B\}$ (odd $p$), long $\{a_1,\dots,a_6\}$, an Association or an equation. \texttt{field}: a prime $p$ or $\{p,r\}$; over $\mathbb F_{p^r}$, $r\ge2$, field elements are written $0$ or $g^k$ (see \texttt{ECPFiniteFieldInfo}).
\returns\ The value $f_{n,P}(Q)\in\mathbb F_q$ of Miller's function, where $\operatorname{div}f_{n,P}=n(P)-([n]P)-(n-1)(\mathcal O)$, normalised as in Miller's algorithm: at every step the line through $T$ and $P$ is divided by the vertical line through $T+P$. An integer for $r=1$, $0$ or $g^k$ for $r\ge2$.
\notes\ Double-and-add loop. Returns \texttt{\$Failed} with a message if $Q$ is a zero or a pole of an intermediate function (for example when $Q$ lies in the group generated by $P$). If $P=\mathcal O$ the result is $1$. Works in every characteristic, including 2 and 3.

\begin{exbox}
\textbf{Example.} $f_{3,P}(Q)$ for $P=(0,4)$, $Q=(5,1)$ on $y^2=x^3+2$ over $\mathbb F_7$

\MatIn{78}{ECPMillerFp[3, \{0, 4\}, \{5, 1\}, \{0, 2\}, 7]}

\smallskip

\MatOut{78}{4}
\end{exbox}

\fname{ECPMinimalModel}
\usageline{ECPMinimalModel[\{a1,...,a6\}]\quad ECPMinimalModel[curve, opts]}
\args\ Any curve form (short $\{A,B\}$ is read as $\{0,0,0,A,B\}$).
\returns\ An Association with \texttt{"MinimalCoefficients"} (normalised so $a_1,a_3\in\{0,1\}$, $a_2\in\{-1,0,1\}$), \texttt{"Equation"}, \texttt{"Transformation"} ($\{u,r,s,t\}$ from the input to the minimal model), \texttt{"Discriminant"}, \texttt{"OldDiscriminant"} and \texttt{"IsMinimal"}.
\notes\ Iterates \texttt{ECPFindLocalReduction} over all primes dividing $\Delta$; feed \texttt{"Transformation"} to \texttt{ECPToMinimalPoint} to carry points along. Its computation may take time, since it factorises $\Delta$, which is slow when $|\Delta|$ has large prime factors.

\begin{exbox}
\textbf{Example.} minimal model, with the transformation used

\MatIn{42}{ECPMinimalModel[\{0, 0, 0, -576, 0\}]}

\smallskip

\MatOut{42}{<|"MinimalCoefficients" -> \{0, 0, 0, -36, 0\}, "Equation" -> y\textasciicircum{}2 == -36 x + x\textasciicircum{}3, "Transformation" -> \{2, 0, 0, 0\}, "Discriminant" -> 2985984, "OldDiscriminant" -> 12230590464, "IsMinimal" -> True|>}
\end{exbox}

\fname{ECPMinimalModelShort}
\usageline{ECPMinimalModelShort[\{A,B\}]\quad ECPMinimalModelShort[curve]}
\args\ A short pair, or any curve form.
\returns\ For short input: the minimal short model (removes $p^4\mid A$, $p^6\mid B$). For long input: the result of \texttt{ECPMinimalModel}.
\notes\ The shape of the Association depends on the kind of input; test for \texttt{"MinimalCoefficients"} when writing generic code.

\begin{exbox}
\textbf{Example.} minimal short model

\MatIn{43}{ECPMinimalModelShort[\{-576, 0\}]}

\smallskip

\MatOut{43}{<|"MinimalCoefficients" -> \{-36, 0\}, "Equation" -> y\textasciicircum{}2 == -36 x + x\textasciicircum{}3, "TransformationScale" -> 2, "Discriminant" -> 2985984, "OldDiscriminant" -> 12230590464, "IsMinimal" -> True|>}
\end{exbox}

\fname{ECPMinimalQ}
\usageline{ECPMinimalQ[\{A,B\}, p]\quad ECPMinimalQ[\{A,B\}]\quad ECPMinimalQ[\{a1,...,a6\}, p]\quad ECPMinimalQ[\{a1,...,a6\}]}
\args\ A short pair or a long model, optionally a prime $p$.
\returns\ \texttt{True}/\texttt{False}. Short overload: not both $p^4\mid A$ and $p^6\mid B$. Long overload: true Laska--Kraus minimality.
\notes\ Short-model minimality is weaker than Weierstrass minimality at $p=2,3$: for instance $y^2=x^3+16$ passes the short test but reduces to $y^2+y=x^3$. Use \texttt{ECPIsMinimal} for the true test.

\begin{exbox}
\textbf{Example.} local test at the prime 2

\MatIn{44}{ECPMinimalQ[\{-576, 0\}, 2]}

\smallskip

\MatOut{44}{False}
\end{exbox}
\medskip\noindent\textbf{Remark (choosing a minimality test).}
\texttt{ECPIsMinimal}, \texttt{ECPIsMinimalShort} and \texttt{ECPMinimalQ} are not interchangeable.
\texttt{ECPIsMinimal} always applies the Laska--Kraus test and answers whether the model is a global minimal Weierstrass model; it takes no prime and a short pair $\{A,B\}$ is read as $\{0,0,0,A,B\}$.
\texttt{ECPMinimalQ} can test at a single prime $p$ (or globally if $p$ is omitted). On a long model it agrees with \texttt{ECPIsMinimal}; on a short pair it only checks that $p^4\mid A$ and $p^6\mid B$ do not both hold, which is minimality among short models.
\texttt{ECPIsMinimalShort} uses this short-model rule for a short pair and equals \texttt{ECPIsMinimal} for long input.
The two notions differ at $p=2,3$: $y^2=x^3+16$ passes the short test, yet it is not minimal as a Weierstrass model (it reduces to $y^2+y=x^3$), so \texttt{ECPIsMinimal[\{0, 16\}]} is \texttt{False} while \texttt{ECPMinimalQ[\{0, 16\}]} is \texttt{True}.
Use \texttt{ECPIsMinimal} whenever the true answer matters, and \texttt{ECPMinimalQ[curve, p]} for a local check at one prime, preferably on a long model.

\fname{ECPMultiply}
\usageline{ECPMultiply[n, P, curve]}
\args\ \texttt{n}: an integer (any sign; $0$ gives \texttt{Infinity}). \texttt{P}: a point or \texttt{Infinity}. \texttt{curve}: any curve form.
\returns\ The point $nP$ by double-and-add, in the coordinates of \texttt{curve}.
\notes\ Negative $n$ uses $-P=(x,-y-a_1x-a_3)$.

\begin{exbox}
\textbf{Example.} 4 (0,0) on 37a1

\MatIn{45}{ECPMultiply[4, \{0, 0\}, \{0, 0, 1, -1, 0\}]}

\smallskip

\MatOut{45}{\{2, -3\}}
\end{exbox}

\fname{ECPMultiplyFp}
\usageline{ECPMultiplyFp[n, P, a, p]\quad ECPMultiplyFp[n, P, curve, p]\quad ECPMultiplyFp[n, P, curve, \{p, r\}]}
\args\ \texttt{n}: an integer of either sign. \texttt{P}, \texttt{a}/\texttt{curve}, \texttt{p}: as in \texttt{ECPAddFp}.
\returns\ $nP$ as $\{x,y\}$ or \texttt{Infinity}.
\notes\ Double-and-add on top of \texttt{ECPAddFp}.

\begin{exbox}
\textbf{Example.} 3 (5,1)

\MatIn{46}{ECPMultiplyFp[3, \{5, 1\}, 2, 17]}

\smallskip

\MatOut{46}{\{10, 6\}}
\end{exbox}

\fname{ECPNaiveHeight}
\usageline{ECPNaiveHeight[\{x, y\}]}
\args\ A rational point (only $x$ is used).
\returns\ $\log\max(|\mathrm{num}(x)|,|\mathrm{den}(x)|,1)$.
\notes\ The base case of the limit in \texttt{ECPCanonicalHeight}; not itself a group invariant. A symbolic $x$ raises a message.

\begin{exbox}
\textbf{Example.} $\log\max(|\mathrm{num}(x)|,|\mathrm{den}(x)|)$ for $x=-3$

\MatIn{47}{ECPNaiveHeight[\{-3, 9\}]}

\smallskip

\MatOut{47}{1.09861}
\end{exbox}

\fname{ECPNormalFormCharTwoThree}
\usageline{ECPNormalFormCharTwoThree[curve, p]\quad ECPNormalFormCharTwoThree[curve, \{p, r\}]}
\args\ \texttt{curve}: any curve form. \texttt{p}: 2 or 3; \texttt{\{p,r\}}: work over $\mathbb F_{p^r}$.
\returns\ An Association with \texttt{"NormalForm"}, \texttt{"Transformation"}, \texttt{"Case"} (a string naming the ordinary or supersingular normal form used) and \texttt{"Characteristic"}. For $p=2$: $y^2+xy=x^3+a_2x^2+a_6$ or $y^2+a_3y=x^3+a_4x+a_6$; for $p=3$: $y^2=x^3+a_2x^2+a_6$ or $y^2=x^3+a_4x+a_6$.
\notes\ Over $\mathbb F_{p^r}$ the result also has \texttt{"FieldOrder"} and, for $p=2$ ordinary, a \texttt{"Gamma"} of absolute trace 1 ($a_2\in\{0,\gamma\}$).

\begin{exbox}
\textbf{Example.} characteristic 2

\MatIn{48}{ECPNormalFormCharTwoThree[\{1, 0, 1, 1, 1\}, 2]}

\smallskip

\MatOut{48}{<|"NormalForm" -> \{1, 1, 0, 0, 1\}, "Transformation" -> \{1, 1, 0, 0\}, "Case" -> "ordinary (a1 != 0): y\textasciicircum{}2+xy=x\textasciicircum{}3+a2 x\textasciicircum{}2+a6", "Characteristic" -> 2|>}
\end{exbox}

\fname{ECPOnCurveQ}
\usageline{ECPOnCurveQ[P, curve]}
\args\ \texttt{P}: a point or \texttt{Infinity}; \texttt{curve}: any curve form.
\returns\ \texttt{True} iff $P$ lies exactly on \texttt{curve}; \texttt{Infinity} always does.
\notes\ Takes a single point list and any curve form, unlike \texttt{ECPPointQ}.

\begin{exbox}
\textbf{Example.} exact membership test, long form

\MatIn{49}{ECPOnCurveQ[\{0, 0\}, \{0, 0, 1, -1, 0\}]}

\smallskip

\MatOut{49}{True}
\end{exbox}

\fname{ECPParseAffineCurve}
\usageline{ECPParseAffineCurve[input]}
\args\ A polynomial, equation, 5-element coefficient list, Association or \texttt{Automatic}.
\returns\ The list $\{F,\{x_v,y_v\}\}$.
\notes\ Internal workhorse of the curve-input handling.

\begin{exbox}
\textbf{Example.} normalise a curve to {polynomial, variables}

\MatIn{50}{ECPParseAffineCurve[y\textasciicircum{}2 + y - x\textasciicircum{}3 + x]}

\smallskip

\MatOut{50}{\{x - x\textasciicircum{}3 + y + y\textasciicircum{}2, \{x, y\}\}}
\end{exbox}

\fname{ECPParseHomogeneousCurve}
\usageline{ECPParseHomogeneousCurve[input]}
\args\ As \texttt{ECPParseAffineCurve}; a two-variable input is homogenised automatically.
\returns\ A homogeneous cubic and its three variables.
\notes\ Used by \texttt{ECPChordAndTangent} and \texttt{ECPCubicToWeierstrass}.

\begin{exbox}
\textbf{Example.} homogenise a cubic with respect to z

\MatIn{51}{ECPParseHomogeneousCurve[y\textasciicircum{}2 + y == x\textasciicircum{}3 - x]}

\smallskip

\MatOut{51}{-x\textasciicircum{}3 + y\textasciicircum{}2 z + x z\textasciicircum{}2 + y z\textasciicircum{}2}
\end{exbox}

\fname{ECPPointOrder}
\usageline{ECPPointOrder[point, curve]}
\args\ A point or \texttt{Infinity}; any curve form.
\returns\ The exact order, or \texttt{Infinity} for a point of infinite order.
\notes\ Over $\mathbb Q$ no torsion order exceeds Mazur's bound 12, so \texttt{Infinity} is unambiguous.

\begin{exbox}
\textbf{Example.} order of (-3,9): infinite, it is not torsion

\MatIn{52}{ECPPointOrder[\{-3, 9\}, \{-36, 0\}]}

\smallskip

\MatOut{52}{Infinity}
\end{exbox}

\fname{ECPPointOrderFp}
\usageline{ECPPointOrderFp[point, a, b, p]\quad ECPPointOrderFp[point, curve, p]\quad ECPPointOrderFp[point, curve, \{p, r\}]}
\args\ A point or \texttt{Infinity} on the curve mod $p$; \texttt{a,b} or \texttt{curve}; a prime $p$ or $\{p,r\}$.
\returns\ The smallest divisor of $N=\#E(\mathbb F_q)$ that kills the point; \texttt{\$Failed} with a message for a non-prime, a singular curve or a point off the curve.
\notes\ Exact, not limited by Mazur's bound. Limited by the point count: about $p\le10^6$, or $q\le$ \texttt{\$ECPFqTableLimit}. Its computation may take time for large $p$, since each call first counts the whole group.

\begin{exbox}
\textbf{Example.} order of (5,1)

\MatIn{53}{ECPPointOrderFp[\{5, 1\}, 2, 2, 17]}

\smallskip

\MatOut{53}{19}
\end{exbox}

\fname{ECPPointQ}
\usageline{ECPPointQ[F, x1, y1]}
\args\ \texttt{F}: a curve Association, equation or polynomial in $x,y$ or $X,Y$. \texttt{x1}, \texttt{y1}: the two coordinates.
\returns\ \texttt{True} iff $(x_1,y_1)$ satisfies \texttt{F}.
\notes\ Compare \texttt{ECPOnCurveQ}, which takes a point list and any curve form.

\begin{exbox}
\textbf{Example.} (0,0) lies on y\textasciicircum{}2 + y = x\textasciicircum{}3 - x

\MatIn{54}{ECPPointQ[y\textasciicircum{}2 + y == x\textasciicircum{}3 - x, 0, 0]}

\smallskip

\MatOut{54}{True}
\end{exbox}

\fname{ECPPointQFp}
\usageline{ECPPointQFp[P, a, b, p]\quad ECPPointQFp[P, curve, p]\quad ECPPointQFp[P, curve, \{p, r\}]}
\args\ A point (not necessarily reduced) or \texttt{Infinity}; \texttt{a,b} or \texttt{curve}; a prime $p$ or $\{p,r\}$.
\returns\ \texttt{True} iff $P$ lies on the curve over the field; \texttt{\$Failed} with a message if $p$ is not prime.
\notes\ Over an extension field give coordinates and coefficients in $0/g^k$ notation.

\begin{exbox}
\textbf{Example.} (5,1) on y\textasciicircum{}2 = x\textasciicircum{}3 + 2x + 2 over F\_17

\MatIn{55}{ECPPointQFp[\{5, 1\}, 2, 2, 17]}

\smallskip

\MatOut{55}{True}
\end{exbox}

\fname{ECPPolynomialHomogeneous}
\usageline{ECPPolynomialHomogeneous[F, vars]}
\args\ A polynomial and its variables.
\returns\ \texttt{True} iff \texttt{F} is homogeneous in \texttt{vars}.
\notes\ No numerical tolerance.

\begin{exbox}
\textbf{Example.} homogeneity test

\MatIn{56}{ECPPolynomialHomogeneous[x\textasciicircum{}2 + x y + y\textasciicircum{}2, \{x, y\}]}

\smallskip

\MatOut{56}{True}
\end{exbox}

\fname{ECPPolynomialTransform}
\usageline{ECPPolynomialTransform[F, vars, matrix]}
\args\ A polynomial, its variables and a square matrix.
\returns\ The expanded polynomial after substituting the linear change of variables.
\notes\ Used when moving a flex to $[0:1:0]$.

\begin{exbox}
\textbf{Example.} linear change of variables

\MatIn{57}{ECPPolynomialTransform[x\textasciicircum{}2 + y\textasciicircum{}2, \{x, y\}, \{\{1, 1\}, \{0, 1\}\}]}

\smallskip

\MatOut{57}{x\textasciicircum{}2 + 2 x y + 2 y\textasciicircum{}2}
\end{exbox}

\fname{ECPPrimeDivisorCount}
\usageline{ECPPrimeDivisorCount[n]\quad ECPPrimeDivisorCount[curve]}
\args\ A nonzero integer $n$, or an integral model in any curve form.
\returns\ $\nu(n)$, the number of distinct primes dividing $n$ ($\nu(\pm1)=0$); for a curve, $\nu(\Delta)$.
\notes\ For a minimal model, $\nu(\Delta)$ is the number of primes of bad reduction. This is the function $\nu$ of Theorem~\ref{thm:rankbound}.

\begin{exbox}
\textbf{Example.} $\nu(360)=3$, and $\nu(\Delta)$ for $y^2=x^3-36x$ ($\Delta=2^{12}3^6$)

\MatIn{99}{\{ECPPrimeDivisorCount[360], ECPPrimeDivisorCount[\{-36, 0\}]\}}

\smallskip

\MatOut{99}{\{3, 2\}}
\end{exbox}

\fname{ECPQuarticToWeierstrass}
\usageline{ECPQuarticToWeierstrass[\{a,b,c,d,e\}]}
\args\ The coefficients of $y^2=ax^4+bx^3+cx^2+dx+e$.
\returns\ A pair $\{A,B\}=\{-27I,-27J\}$ with $I=12ae-3bd+c^2$, $J=72ace+9bcd-27ad^2-27eb^2-2c^3$: the Jacobian of the quartic.
\notes\ The Jacobian is isomorphic to the quartic itself exactly when the quartic has a rational point.

\begin{exbox}
\textbf{Example.} y\textasciicircum{}2 = x\textasciicircum{}4 + x + 1

\MatIn{58}{ECPQuarticToWeierstrass[\{1, 0, 0, 1, 1\}]}

\smallskip

\MatOut{58}{\{-324, 729\}}
\end{exbox}

\fname{ECPRankBound}
\usageline{ECPRankBound[curve]\quad ECPRankBound[curve, "Details"]}
\args\ An integral model in any curve form with a rational point of order 2.
\returns\ The smallest of the proven bounds of Theorem~\ref{thm:rankbound} and Proposition~\ref{prop:threeroots} that apply: $\nu(a^2-4b)+\nu(b)-1$ after moving a rational root of the cubic to $x=0$, $m+2\alpha-1$ from the bad reduction, and $2\nu(\Delta)$ when all three 2-torsion points are rational. With \texttt{"Details"}, an Association with \texttt{"Bound"}, \texttt{"Candidates"} (method $\to$ bound) and the data used.
\notes\ A proven upper bound, not the rank. If $E(\mathbb Q)$ has no point of order 2 no bound of this type exists, and \texttt{Missing["NoRational2Torsion"]} is returned with a message. Its computation may take time, since \texttt{ECPLocalData} factorises $\Delta$.

\begin{exbox}
\textbf{Example.} $y^2=x^3-36x$ (rank 1), and 37a1, which has no rational point of order 2

\MatIn{100}{\{ECPRankBound[\{-36, 0\}], ECPRankBound[\{0, 0, 1, -1, 0\}]\}}

\smallskip

\MatOut{100}{\{3, Missing["NoRational2Torsion"]\}}
\end{exbox}

\fname{ECPRankUpperBound}
\usageline{ECPRankUpperBound[\{A,B\}]\quad ECPRankUpperBound[curve]}
\args\ Any curve form with integral coefficients.
\returns\ A small integer: the number of distinct prime factors of $\Delta$, capped at 4.
\notes\ A crude sanity signal, not derived from 2-descent: it equals the rank for 37a1 but overestimates it for $y^2=x^3-x$. Its computation may take time, since it factorises $\Delta$, which is slow when $|\Delta|$ has large prime factors.

\begin{exbox}
\textbf{Example.} heuristic rank estimate (not a proven bound)

\MatIn{59}{ECPRankUpperBound[\{-36, 0\}]}

\smallskip

\MatOut{59}{2}
\end{exbox}

\fname{ECPSchoofPointCount}
\usageline{ECPSchoofPointCount[a, b, p]\quad ECPSchoofPointCount[curve, p]}
\args\ \texttt{a,b} or \texttt{curve}: short $\{A,B\}$, long $\{a_1,\dots,a_6\}$, an Association or an equation; a prime $p$ of any size ($\{p,1\}$ is also accepted).
\returns\ $\#E(\mathbb F_p)$, including the point at infinity.
\notes\ Schoof's algorithm: $a_p \bmod \ell$ is found from the action of the Frobenius endomorphism on $E[\ell]$ (division polynomials, arithmetic in $\mathbb F_p[x]/(\psi_\ell)$) for $\ell=2,3,5,\dots$ ($\ell\neq p$) until $\prod\ell>4\sqrt p$, and $a_p$ is recovered by the Chinese Remainder Theorem and Hasse's bound. The running time is polynomial in $\log p$, so $p$ is not limited to $10^6$ as for \texttt{ECPSEAPointCount}. For $p=2,3$ the points are counted directly. The Elkies and Atkin improvements (SEA) are not implemented. Fields $\mathbb F_{p^r}$ with $r>1$ are not supported; use \texttt{ECPSEAPointCount}. Its computation may take time for very large $p$, since the division polynomials $\psi_\ell$ have degree about $\ell^2/2$.

\begin{exbox}
\textbf{Example.} a prime far beyond the enumeration limit, and a check against enumeration

\MatIn{83}{ECPSchoofPointCount[\{2, 3\}, 1000000000039]}

\smallskip

\MatOut{83}{1000000769568}

\medskip

\MatIn{84}{ECPSchoofPointCount[\{2, 3\}, 1000003] == ECPSEAPointCount[\{2, 3\}, 1000003]}

\smallskip

\MatOut{84}{True}
\end{exbox}

\fname{ECPSEAPointCount}
\usageline{ECPSEAPointCount[a, b, p]\quad ECPSEAPointCount[curve, p]\quad ECPSEAPointCount[curve, \{p, r\}]}
\args\ \texttt{a,b} or \texttt{curve} (the long overload also works for $p=2,3$); a prime $p\lesssim10^6$ or $\{p,r\}$.
\returns\ $\#E(\mathbb F_q)$ including the point at infinity.
\notes\ Despite the name this is brute-force enumeration, not the Schoof--Elkies--Atkin algorithm; for larger primes use \texttt{ECPSchoofPointCount}. For a curve with plain integer coefficients over $\mathbb F_{p^r}$ the count comes from $\#E(\mathbb F_p)$ by a Lucas-sequence recurrence and works for any $r$. Its computation may take time for large $p$, since the running time is proportional to $p$.

\begin{exbox}
\textbf{Example.} number of points (direct count, not the SEA algorithm)

\MatIn{60}{ECPSEAPointCount[2, 2, 17]}

\smallskip

\MatOut{60}{19}
\end{exbox}

\fname{ECPSearchPoints}
\usageline{ECPSearchPoints[curve, bound]}
\args\ \texttt{curve}: any curve form; \texttt{bound} (default 50): numerators $|x|\le$ \texttt{bound}, denominators up to 10.
\returns\ All points found, in the coordinates of \texttt{curve}.
\notes\ Not guaranteed complete for positive rank; method distinct from Nagell--Lutz.

\begin{exbox}
\textbf{Example.} rational points of small height

\MatIn{61}{ECPSearchPoints[\{-36, 0\}, 30]}

\smallskip

\MatOut{61}{\{\{-6, 0\}, \{-3, 9\}, \{-3, -9\}, \{-2, 8\}, \{-2, -8\}, \{0, 0\}, \{6, 0\}, \{12, 36\}, \{12, -36\}, \{18, 72\}, \{18, -72\}, \{25/4, 35/8\}, \{25/4, -35/8\}\}}
\end{exbox}

\fname{ECPSetEC}
\usageline{ECPSetEC[a1, a2, a3, a4, a6]}
\args\ The five long coefficients as separate arguments.
\returns\ The new value of \texttt{\$ECPCurrentEC}.
\notes\ A singular input raises a message and leaves \texttt{\$ECPCurrentEC} unchanged.

\begin{exbox}
\textbf{Example.} set the session curve to 37a1: y\textasciicircum{}2 + y = x\textasciicircum{}3 - x

\MatIn{62}{ECPSetEC[0, 0, 1, -1, 0]}

\smallskip

\MatOut{62}{<|"Coefficients" -> \{0, 0, 1, -1, 0\}, "Equation" -> x - x\textasciicircum{}3 + y + y\textasciicircum{}2 == 0|>}
\end{exbox}

\fname{ECPShortDiscriminant}
\usageline{ECPShortDiscriminant[\{A,B\}]\quad ECPShortDiscriminant[A, B]}
\args\ A short pair or two numbers.
\returns\ $-16(4A^3+27B^2)$.
\notes\ Short input only; use \texttt{ECPWeierstrassDiscriminant} for other forms.

\begin{exbox}
\textbf{Example.} discriminant of y\textasciicircum{}2 = x\textasciicircum{}3 - 36x

\MatIn{63}{ECPShortDiscriminant[\{-36, 0\}]}

\smallskip

\MatOut{63}{2985984}
\end{exbox}

\fname{ECPShortWeierstrassForm}
\usageline{ECPShortWeierstrassForm[curve]}
\args\ Any curve form.
\returns\ Curve in Short form with integer coefficients.
\notes\ Uses a different transformation than \texttt{ECPToShortWeierstrass}.

\begin{exbox}
\textbf{Example.} correctly spelled synonym

\MatIn{64}{ECPShortWeierstrassForm[y\textasciicircum{}2 + y == x\textasciicircum{}3 - x]}

\smallskip

\MatOut{64}{-11664 + 1296 x - x\textasciicircum{}3 + y\textasciicircum{}2}
\end{exbox}

\fname{ECPTate}
\usageline{ECPTate[model, p]}
\args\ An integral long model $\{a_1,\dots,a_6\}$ and a prime $p$ (normally a bad prime).
\returns\ An Association with the conductor exponent (key \texttt{"f"}), the Kodaira symbol, the minimal model at $p$ and the split type.
\notes\ At a prime of good reduction it only confirms minimality and returns exponent 0. A singular model or a non-prime raises a message.

\begin{exbox}
\textbf{Example.} Tate's algorithm for 37a1 at p = 37

\MatIn{66}{ECPTate[\{0, 0, 1, -1, 0\}, 37]}

\smallskip

\MatOut{66}{<|"f" -> 1, "type" -> "I1", "model" -> \{0, 15, 37, 74, -222\}, "split" -> False|>}
\end{exbox}

\fname{ECPTateForm}
\usageline{ECPTateForm[orders, t]}
\args\ \texttt{orders}: a list of five non-negative integers $(o_1,\dots,o_5)$ for $(a_1,a_2,a_3,a_4,a_6)$. \texttt{t}: the parameter (a symbol).
\returns\ The long Weierstrass coefficients $a_i=(c_i+(i+1)\,t)\,t^{o_i}$ with the base constants $c=(1,1,3,5,7)$, a model whose coefficients vanish to the prescribed orders at $t=0$.
\notes\ Example orders: SU(2) $(0,0,1,1,2)$, SU(3) $(0,1,1,2,3)$, SU(5) $(0,1,2,3,5)$, SO(10) $(1,1,2,3,5)$, E6 $(1,2,2,3,5)$, E7 $(1,2,3,3,5)$, E8 $(1,2,3,4,5)$. The base constants are one arbitrary generic choice, not the most general Tate form.

\begin{exbox}
\textbf{Example.} Tate form for $\mathfrak{su}(2)$

\MatIn{96}{ECPTateForm[\{0, 0, 1, 1, 2\}, t]}

\smallskip

\MatOut{96}{\{1 + 2 t, 1 + 3 t, t (3 + 4 t), t (5 + 5 t), t\textasciicircum{}2 (7 + 6 t)\}}
\end{exbox}

\fname{ECPTatePairingFp}
\usageline{ECPTatePairingFp[n, P, Q, curve, field]}
\args\ \texttt{n}: a positive integer dividing $q-1$. \texttt{P}: a point of $E(\mathbb F_q)$ with $nP=\mathcal O$. \texttt{Q}: any point of $E(\mathbb F_q)$. \texttt{curve}: any curve form (every characteristic). \texttt{field}: a prime $p$ or $\{p,r\}$.
\returns\ The reduced Tate pairing $t_n(P,Q)=f_{n,P}(Q)^{(q-1)/n}$, an $n$-th root of unity in $\mathbb F_q$ (an integer for $r=1$, $0$ or $g^k$ for $r\ge2$).
\notes\ Computed by Miller's algorithm. The value depends only on $Q$ modulo $nE(\mathbb F_q)$; the pairing is bilinear and non-degenerate on $E(\mathbb F_q)[n]\times E(\mathbb F_q)/nE(\mathbb F_q)$. A point $P$ of smaller exact order $m$ is reduced to $t_m(P,Q)$. If $Q$ is a zero or a pole of an intermediate Miller function, the value $f(Q+S)/f(S)$ with an auxiliary point $S$ is used. If $n\nmid q-1$ a message asks for an extension field with $n\mid p^r-1$. Unlike the Weil pairing the pairing is not alternating, so $t_n(P,P)$ need not be $1$. Its computation may take time for large $q$ when the auxiliary point $S$ is needed, since the points of $E(\mathbb F_q)$ are then listed.

\begin{exbox}
\textbf{Example.} Tate pairing of $P=(0,4)$ and $Q=(5,1)$ on $y^2=x^3+2$ over $\mathbb F_7$

\MatIn{79}{ECPTatePairingFp[3, \{0, 4\}, \{5, 1\}, \{0, 2\}, 7]}

\smallskip

\MatOut{79}{2}
\end{exbox}

\fname{ECPToLongCoordinates}
\usageline{ECPToLongCoordinates[pt, curve]\quad ECPToLongCoordinates[pt, curve, \{u,r,s,t\}]}
\args\ A point $(x',y')$ of the transformed model \texttt{ECPTransformWeierstrass[curve, \{u,r,s,t\}]}; the original long curve; an isomorphism $\{u,r,s,t\}$ (default identity $\{1,0,0,0\}$). The point is checked on the transformed model.
\returns\ The point on the original model via $x=u^2x'+r$, $y=u^3y'+u^2sx'+t$.
\notes\ With the default identity this is a no-op; to undo \texttt{ECPToShortCoordinates} use \texttt{ECPFromShortCoordinates}.

\begin{exbox}
\textbf{Example.} $(0,0)$ on the transformed model $\{0,3,1,2,0\}$ (see \texttt{ECPTransformWeierstrass}) is $(1,0)$ on 37a1

\MatIn{67}{ECPToLongCoordinates[\{0, 0\}, \{0, 0, 1, -1, 0\}, \{1, 1, 0, 0\}]}

\smallskip

\MatOut{67}{\{1, 0\}}
\end{exbox}

\fname{ECPToMinimalPoint}
\usageline{ECPToMinimalPoint[\{x,y\}, \{u,r,s,t\}]}
\args\ A point of the original model and the isomorphism (e.g.\ \texttt{"Transformation"} of \texttt{ECPMinimalModel}).
\returns\ The point on the new model: $x'=(x-r)/u^2$, $y'=(y-u^2sx'-t)/u^3$.
\notes\ $u=0$ raises a message.

\begin{exbox}
\textbf{Example.} map (-12, 72) on y\textasciicircum{}2 = x\textasciicircum{}3 - 576x to the minimal model

\MatIn{68}{ECPToMinimalPoint[\{-12, 72\}, ECPMinimalModel[\{0, 0, 0, -576, 0\}]["Transformation"]]}

\smallskip

\MatOut{68}{\{-3, 9\}}
\end{exbox}

\fname{ECPTorsionBasisFp}
\usageline{ECPTorsionBasisFp[n, curve, field]}
\args\ \texttt{n}: an integer $\ge2$ prime to $p$. \texttt{curve}: any curve form. \texttt{field}: a prime $p$ or $\{p,r\}$.
\returns\ A basis $\{P,Q\}$ of the full $n$-torsion subgroup $E[n]\subset E(\mathbb F_q)$: $P$ has exact order $n$ and \texttt{ECPWeilPairingFp[n, P, Q, curve, field]} is a primitive $n$-th root of unity.
\notes\ Requires $E[n]\subset E(\mathbb F_q)$, that is, $E(\mathbb F_q)\supseteq\mathbb Z/n\times\mathbb Z/n$ (this forces $n\mid q-1$); otherwise a message says so, suggests an extension field $\{p,r\}$, and \texttt{\$Failed} is returned. All points are enumerated, so $q\le$ \texttt{\$ECPFqTableLimit}. Over $\mathbb F_{p^r}$ the coordinates are written $0$ or $g^k$. Its computation may take time for large $q$, since every point of $E(\mathbb F_q)$ is listed and multiplied by $n$.

\begin{exbox}
\textbf{Example.} basis of $E[3]$ for $y^2=x^3+2$ over $\mathbb F_7$

\MatIn{80}{ECPTorsionBasisFp[3, \{0, 2\}, 7]}

\smallskip

\MatOut{80}{\{\{0, 4\}, \{5, 1\}\}}
\end{exbox}

\fname{ECPTorsionFamily}
\usageline{ECPTorsionFamily[N, d]}
\args\ \texttt{N}: 4, 5, 6 or 7. \texttt{d}: the parameter on $X_1(N)$ (a symbol).
\returns\ The universal elliptic curve with a rational $N$-torsion point, in the Tate--Kubert normal form $y^2+(1-c)xy-by=x^3-bx^2$ (the point $(0,0)$ has order $N$), as an Association with the keys \texttt{"Coefficients"} ($\{1-c,-b,-b,0,0\}$), \texttt{"b"}, \texttt{"c"} and \texttt{"Parameter"}.
\notes\ Over $\mathbb Q(d)$ the torsion of the Mordell--Weil group contains $\mathbb Z/N\mathbb Z$. Use \texttt{ECPFibreList} on the coefficients to see the fibre content ($N=6$ gives $I_6,I_3,I_2,I_1$). Other values of $N$ raise a message.

\begin{exbox}
\textbf{Example.} the family with a point of order 6

\MatIn{97}{ECPTorsionFamily[6, d]["Coefficients"]}

\smallskip

\MatOut{97}{\{1 - d, -d - d\textasciicircum{}2, -d - d\textasciicircum{}2, 0, 0\}}
\end{exbox}

\fname{ECPTorsionSubgroup}
\usageline{ECPTorsionSubgroup[curve]}
\args\ Any curve form.
\returns\ An Association with \texttt{"Order"}, \texttt{"Structure"} (a string), \texttt{"Generators"} and \texttt{"Points"} (including \texttt{Infinity}); points are returned in the coordinates of \texttt{curve}.
\notes\ Built on \texttt{ECPFindTorsionPoints}, so it inherits the \texttt{\$ECPNagellLutzDiscLimit} fallback. Its computation may take time for the same reason as \texttt{ECPFindTorsionPoints}.

\begin{exbox}
\textbf{Example.} torsion subgroup of y\textasciicircum{}2 = x\textasciicircum{}3 - 36x

\MatIn{69}{ECPTorsionSubgroup[\{-36, 0\}]}

\smallskip

\MatOut{69}{<|"Order" -> 4, "Structure" -> "Z/2Z x Z/2Z", "Generators" -> \{\{0, 0\}, \{-6, 0\}\}, "Points" -> \{Infinity, \{-6, 0\}, \{0, 0\}, \{6, 0\}\}|>}
\end{exbox}

\fname{ECPToShortCoordinates}
\usageline{ECPToShortCoordinates[\{x,y\}, curve]}
\args\ A point on \texttt{curve}; any curve form.
\returns\ The point $\{X,Y\}$ on the model of \texttt{ECPToShortWeierstrass[curve]}: $X=x+b_2/12$, $Y=y+(a_1x+a_3)/2$.
\notes\ A point off the curve, or a singular curve, raises a message and returns \texttt{\$Failed}.

\begin{exbox}
\textbf{Example.} move the point (0,0) to the short model

\MatIn{70}{ECPToShortCoordinates[\{0, 0\}, \{0, 0, 1, -1, 0\}]}

\smallskip

\MatOut{70}{\{0, 1/2\}}
\end{exbox}

\fname{ECPToShortCoordinatesFp}
\usageline{ECPToShortCoordinatesFp[\{x,y\}, curve, p]}
\args\ A point of the long curve over $\mathbb F_p$; any curve form; a prime $p$.
\returns\ The image on the model of \texttt{ECPToShortWeierstrassFp[curve, p]}; \texttt{Infinity} maps to \texttt{Infinity}.
\notes\ A group isomorphism. The point must lie on the curve mod $p$.

\begin{exbox}
\textbf{Example.} move (0,0) to the reduced model

\MatIn{71}{ECPToShortCoordinatesFp[\{0, 0\}, \{0, 0, 1, -1, 0\}, 7]}

\smallskip

\MatOut{71}{\{0, 4\}}
\end{exbox}

\fname{ECPToShortWeierstrass}
\usageline{ECPToShortWeierstrass[curve]}
\args\ Any curve form.
\returns\ The short model $Y^2=X^3+AX+B$ with $A=-c_4/48$, $B=-c_6/864$, as an equation in $X,Y$.
\notes\ For $p=2,3$ use \texttt{ECPToShortWeierstrassFp} or \texttt{ECPNormalFormCharTwoThree}.

\begin{exbox}
\textbf{Example.} the short model Y\textasciicircum{}2 = X\textasciicircum{}3 + A X + B, A = -c4/48, B = -c6/864

\MatIn{72}{ECPToShortWeierstrass[\{0, 0, 1, -1, 0\}]}

\smallskip

\MatOut{72}{Y\textasciicircum{}2 == 1/4 - X + X\textasciicircum{}3}
\end{exbox}
\medskip\noindent\textbf{Remark (\texttt{ECPToShortWeierstrass} versus \texttt{ECPShortWeierstrassForm}).}
Both reduce a long model to short Weierstrass form for characteristic $0$ or $>3$, but they return different things.
\texttt{ECPToShortWeierstrass[curve]} gives the equation $Y^2=X^3+AX+B$ with $A=-c_4/48$, $B=-c_6/864$, written in the new variables $X,Y$; the coefficients may be rational, and the substitution $X=x+b_2/12$, $Y=y+(a_1x+a_3)/2$ is the one used by \texttt{ECPToShortCoordinates} and \texttt{ECPFromShortCoordinates}. For 37a1 it returns \texttt{Y\textasciicircum{}2 == 1/4 - X + X\textasciicircum{}3}.
\texttt{ECPShortWeierstrassForm[curve]} returns the polynomial $y^2-(x^3+Ax+B)$ in the curve's own variables $x,y$, not an equation. For 37a1 it returns $-11664+1296x-x^3+y^2$, that is $y^2=x^3-1296x+11664$: the same curve rescaled by $u=6$ so that the coefficients are integers ($A=-27c_4$, $B=-54c_6$).
The two outputs are therefore isomorphic over $\mathbb Q$ but are not the same equation, and points computed on one model cannot be used on the other without the corresponding scaling. Use \texttt{ECPToShortWeierstrass} when you need the point transformation or the standard $A=-c_4/48$ normalisation, and \texttt{ECPShortWeierstrassForm} when you want an integral short model as a polynomial. Neither applies in characteristic $2$ or $3$; use \texttt{ECPToShortWeierstrassFp}.

\fname{ECPToShortWeierstrassFp}
\usageline{ECPToShortWeierstrassFp[curve, p]\quad ECPToShortWeierstrassFp[curve, \{p, r\}]}
\args\ Any curve form (integer or $p$-integral rational coefficients); a prime $p$ or $\{p,r\}$.
\returns\ An Association with \texttt{"Characteristic"}, \texttt{"OriginalCoefficients"}, \texttt{"ShortCoefficients"}, \texttt{"Transformation"}, \texttt{"Case"}, \texttt{"PureShortQ"}, \texttt{"ShortPair"}, \texttt{"Equation"}, \texttt{"Discriminant"} and \texttt{"jInvariant"} (plus \texttt{"FieldOrder"} over $\mathbb F_{p^r}$).
\notes\ For $p>3$ the result is a genuine short model; for $p=3$ one exists only if $j=0$ (otherwise $y^2=x^3+a_2x^2+a_6$); for $p=2$ none exists.

\begin{exbox}
\textbf{Example.} long model of 37a1 reduced modulo 7

\MatIn{73}{ECPToShortWeierstrassFp[\{0, 0, 1, -1, 0\}, 7]}

\smallskip

\MatOut{73}{<|"Characteristic" -> 7, "OriginalCoefficients" -> \{0, 0, 1, 6, 0\}, "ShortCoefficients" -> \{0, 0, 0, 6, 2\}, "Transformation" -> \{1, 0, 0, 3\}, "Case" -> "p > 3: genuine short model y\textasciicircum{}2 = x\textasciicircum{}3 + A x + B (A = -c4/48, B = -c6/864 mod p)", "PureShortQ" -> True, "ShortPair" -> \{6, 2\}, "Equation" -> y\textasciicircum{}2 == 2 + 6 x + x\textasciicircum{}3, "Discriminant" -> 2, "jInvariant" -> 3|>}
\end{exbox}

\fname{ECPToWeierstrassForm}
\usageline{ECPToWeierstrassForm[F, \{xv, yv\}]}
\args\ A polynomial or equation and its two variables.
\returns\ The equation $y^2+a_1xy+a_3y=x^3+a_2x^2+a_4x+a_6$.
\notes\ No change of variables; stray non-Weierstrass monomials are reported, not discarded.

\begin{exbox}
\textbf{Example.} read the coefficients off a cubic already in Weierstrass shape

\MatIn{74}{ECPToWeierstrassForm[y\textasciicircum{}2 + y - x\textasciicircum{}3 + x, \{x, y\}]}

\smallskip

\MatOut{74}{x - x\textasciicircum{}3 + y + y\textasciicircum{}2}
\end{exbox}

\fname{ECPTransformWeierstrass}
\usageline{ECPTransformWeierstrass[\{a1,...,a6\}, \{u,r,s,t\}]}
\args\ A long model (or any curve form) and an isomorphism with $u\neq0$.
\returns\ The new coefficients $\{a_1',\dots,a_6'\}$ for $x=u^2x'+r$, $y=u^3y'+u^2sx'+t$:
\begin{align*}
  ua_1'&=a_1+2s, \qquad u^2a_2'=a_2-sa_1+3r-s^2, \qquad u^3a_3'=a_3+ra_1+2t,\\
  u^4a_4'&=a_4-sa_3+2ra_2-(t+rs)a_1+3r^2-2st,\\
  u^6a_6'&=a_6+ra_4+r^2a_2+r^3-ta_3-t^2-rta_1.
\end{align*}
\notes\ The building block of \texttt{ECPMinimalModel}, \texttt{ECPFindLocalReduction}, \texttt{ECPNormalFormCharTwoThree} and \texttt{ECPToShortWeierstrassFp}. $u=0$ raises a message.

\begin{exbox}
\textbf{Example.} apply the isomorphism {1,1,0,0} to the coefficients

\MatIn{75}{ECPTransformWeierstrass[\{0, 0, 1, -1, 0\}, \{1, 1, 0, 0\}]}

\smallskip

\MatOut{75}{\{0, 3, 1, 2, 0\}}
\end{exbox}

\fname{ECPWeierstrassDiscriminant}
\usageline{ECPWeierstrassDiscriminant[curve]}
\args\ Any curve form.
\returns\ $\Delta=-b_2^2b_8-8b_4^3-27b_6^2+9b_2b_4b_6$.
\notes\ The curve is nonsingular iff $\Delta\neq0$.

\begin{exbox}
\textbf{Example.} discriminant of 37a1

\MatIn{76}{ECPWeierstrassDiscriminant[\{0, 0, 1, -1, 0\}]}

\smallskip

\MatOut{76}{37}
\end{exbox}

\fname{ECPWeilPairingFp}
\usageline{ECPWeilPairingFp[n, P, Q, curve, field]}
\args\ \texttt{n}: a positive integer prime to $p$. \texttt{P}, \texttt{Q}: points of $E[n]$ ($nP=nQ=\mathcal O$) on \texttt{curve} over the field. \texttt{curve}: any curve form (every characteristic). \texttt{field}: a prime $p$ or $\{p,r\}$.
\returns\ The Weil pairing $e_n(P,Q)=(-1)^n\,f_{n,P}(Q)/f_{n,Q}(P)$, computed with Miller's algorithm for both functions: an $n$-th root of unity in $\mathbb F_q$, written as an integer for the field $p$ or $\{p,1\}$ and as $0$ or $g^k$ for $r\ge2$.
\notes\ The pairing is bilinear, alternating ($e_n(P,P)=1$) and non-degenerate. Points of smaller exact order are reduced to their exact order, and dependent points give $1$. Returns \texttt{\$Failed} with a message if $n$ is not prime to $p$, or if a point is not on the curve or not in $E[n]$.

\begin{exbox}
\textbf{Example.} Weil pairing of the basis of $E[3]$ for $y^2=x^3+2$ over $\mathbb F_7$ (a primitive cube root of unity: $2^3\equiv1$)

\MatIn{81}{ECPWeilPairingFp[3, \{0, 4\}, \{5, 1\}, \{0, 2\}, 7]}

\smallskip

\MatOut{81}{2}
\end{exbox}

\fname{ECPWInvariants}
\usageline{ECPWInvariants[curve]}
\args\ Any curve form.
\returns\ The Association \texttt{<|"b2"->..., "b4"->..., "b6"->..., "b8"->..., "c4"->..., "c6"->..., "Discriminant"->...|>}.
\notes\ There is no $j$ key; use \texttt{ECPJInvariant}.

\begin{exbox}
\textbf{Example.} b2, b4, b6, b8, c4, c6 and the discriminant

\MatIn{77}{ECPWInvariants[\{0, 0, 1, -1, 0\}]}

\smallskip

\MatOut{77}{<|"b2" -> 0, "b4" -> -2, "b6" -> 1, "b8" -> -1, "c4" -> 48, "c6" -> -216, "Discriminant" -> 37|>}
\end{exbox}

\section{Notebook guide}\label{app: nb guide}
Here are some pre-computed example notebooks from various standard textbooks on the elliptic curves using ECPack.wl for the better understanding of the subject. User can click on the name of example notebook, which will redirect them to Github page. There the notebook can be downloaded to check the computed examples.
\begin{table}[htbp]
\centering
\renewcommand{\arraystretch}{1.2}
\begin{tabular}{|p{0.28\textwidth}|p{0.64\textwidth}|}
\hline
\textbf{Notebook} & \textbf{Description} \\
\hline
\href{https://github.com/DhamakedarDarsh/Mathematica-Code-and-worked-out-notebooks-for-Elliptic-Curves/blob/main/Lozano_Examples.nb}{Lozano Examples}~\cite{LozanoRobledo2011} & Textbook-style exercises illustrating changes of form, invariants, and basic computations with elliptic curves in the style of Lozano-Robledo. \\
\hline
\href{https://github.com/DhamakedarDarsh/Mathematica-Code-and-worked-out-notebooks-for-Elliptic-Curves/blob/main/Husemoller_Examples.nb}{Husem{\"o}ller Examples}~\cite{Husemoller2004EllipticCurves} & Classical exercises on the group law, torsion points, and explicit arithmetic drawn from the textbook of Husem{\"o}ller. \\
\hline
\href{https://github.com/DhamakedarDarsh/Mathematica-Code-and-worked-out-notebooks-for-Elliptic-Curves/blob/main/Washington_Examples.nb}{Washington Examples}~\cite{Washington2008EllipticCurves} & Worked examples emphasizing finite-field arithmetic, point counts, and cryptography-related computations from the textbook of Washington. \\
\hline
\href{https://github.com/DhamakedarDarsh/Mathematica-Code-and-worked-out-notebooks-for-Elliptic-Curves/blob/main/Silverman_AEC_Examples.nb}{Silverman AEC Examples}~\cite{Silverman2009ArithmeticEllipticCurves} & Standard arithmetic examples covering minimal models, heights, ranks, and rational points from the textbook of Silverman. \\
\hline
\href{https://github.com/DhamakedarDarsh/Mathematica-Code-and-worked-out-notebooks-for-Elliptic-Curves/blob/main/SilvermanTate_Examples.nb}{Silverman Tate Examples}~\cite{silverman} & Examples of point arithmetic and other computations over rational and finite fields from the textbook of Silverman and Tate. \\
\hline
\href{https://github.com/DhamakedarDarsh/Mathematica-Code-and-worked-out-notebooks-for-Elliptic-Curves/blob/main/McGee_Examples.nb}{McGee Examples}~\cite{McGee2015} & Comparison examples inspired by McGee's package, highlighting point arithmetic over $\mathbb F_p$, Bachet's duplication formula and ElGamal encryption. \\
\hline
\href{https://github.com/DhamakedarDarsh/Mathematica-Code-and-worked-out-notebooks-for-Elliptic-Curves/blob/main/Miscellaneous_Examples.nb}{Miscellaneous Examples} & Additional test curves and miscellaneous examples. Covers all possible commands of ECPack on some chosen curves. \\
\hline
\end{tabular}
\caption{Summary of the seven accompanying Mathematica notebooks.}
\label{tab:notebooks}
\end{table}

\clearpage
\section*{Declarations}

\subsection*{Funding}
This research did not receive any specific funding.

\subsection*{Availability of data and material}
Not applicable. This manuscript does not report data generation or analysis. The package \texttt{ECPack.wl} and the accompanying Mathematica notebooks are openly available in the \href{https://github.com/DhamakedarDarsh/Mathematica-Code-and-worked-out-notebooks-for-Elliptic-Curves}{ECPack GitHub repository}.

\subsection*{Competing interests}
The authors have no competing interests to declare that are relevant to the content of this article.

\subsection*{Authors' contributions}
The Mathematica package \texttt{ECPack.wl} was created by D.B., with help from P.T. in the later stages. B.A. provided the essential knowledge of the subject without which this work would not have been completed, and gave D.B. and P.T. the guidance and mentorship they needed. The manuscript was written by B.A., D.B. and P.T.

\printbibliography

\end{document}